\documentclass[11pt]{article}
\usepackage{ebezier}
\usepackage{amssymb,amsmath,mathrsfs,color,graphicx,float,ifpdf} 
 \allowdisplaybreaks

\catcode`!=11
\let\!int\int \def\int{\displaystyle\!int}
\let\!lim\lim \def\lim{\displaystyle\!lim}
\let\!sum\sum \def\sum{\displaystyle\!sum}
\let\!sup\sup \def\sup{\displaystyle\!sup}
\let\!inf\inf \def\inf{\displaystyle\!inf}
\let\!cap\cap \def\cap{\displaystyle\!cap}
\let\!cup\cup \def\cup{\displaystyle\!cup}
\let\!max\max \def\max{\displaystyle\!max}
\let\!min\min \def\min{\displaystyle\!min}
\let\!frac\frac \def\frac{\displaystyle\!frac}
\catcode`!=12

\newtheorem{theorem}{\bf Theorem}[section]
\newtheorem{lemma}{\bf Lemma}[section]
\newtheorem{definition}{\bf Definition}[section]
\newtheorem{proposition}{\bf Proposition}[section]
\newtheorem{remark}{\bf Remark}[section]

\def\Proof{\it{Proof.}\rm\quad}

\def\pd#1#2{\frac{\partial#1}{\partial#2}}

\catcode`@=11
\@addtoreset{equation}{section}
\catcode`@=12

\begin{document}

\title{On an Elliptic Free Boundary Problem and Subsonic Jet Flows for a Given Surrounding Pressure}
\author{Chunpeng Wang
\footnote{Supported by a grant
from the National Natural Science Foundation of China (No. 11571137).}
\quad \small (email: wangcp@jlu.edu.cn)
\\
\small School of Mathematics, Jilin University, Changchun
130012, China
\\[2mm]
Zhouping Xin\footnote{Supported by
Zheng Ge Ru Foundation, Hong Kong RGC Earmarked Research
Grants CUHK-14305315, CUHK-14300917 and CUHK-14302917,
NSFC/RGC Joint Research Scheme Grant N-CUHK443/14,
and a Focus Area Grant from the Chinese University of Hong Kong.}
\quad
\small (email: zpxin@ims.cuhk.edu.hk)
\\
\small The Institute of Mathematical Sciences and Department of
Mathematics,
\\
\small The Chinese University of Hong Kong, Shatin, NT, Hong Kong
}
\date{}
\maketitle

\begin{abstract}
This paper concerns compressible subsonic jet flows for a given surrounding pressure
from a two-dimensional finitely long convergent nozzle with straight solid wall,
which are governed by a free boundary problem for a quasilinear elliptic equation.
For a given surrounding pressure and a given incoming mass flux,
we seek a subsonic jet flow with the given incoming mass flux such that
the flow velocity at the inlet is along the normal direction,
the flow satisfies the slip condition at the wall,
and the pressure of the flow at the free boundary coincides with the given surrounding pressure.
In general, the free boundary contains two parts: one is the particle path connected with the wall
and the other is a level set of the velocity potential.
We identify a suitable space of flows in terms of
the minimal speed and the maximal velocity potential difference
for the well-posedness of the problem.
It is shown that there is an optimal interval
such that there exists a unique subsonic jet flow in the space
iff the length of the nozzle belongs to this interval.
Furthermore, the optimal regularity and other properties of the flows are shown.
\\[6pt]
{\sl Keywords:} Free boundary, Jet flow, Mixed boundary conditions.
\\
{\sl 2010 MR Subject Classification:} 35R35 76N10 35J66
\end{abstract}

\section{Introduction}

In this paper we study the compressible subsonic jet flows for a given surrounding pressure.
Such problems arise naturally in physical experiments and engineering
designs (\cite{Bers}), and have received much attention for a long time.
Many examples and numerical results can be found
in the monographs \cite{Bers,book1,Chaplygin,book2}.
In general, compressible subsonic jet flows are governed by elliptic free boundary problems.
The first rigorous mathematical theory was established until 1980's.
H. W. Alt, L. A. Caffarelli and A. Friedman developed a variational approach to solve
free boundary problem for elliptic equations in \cite{Alt1,Alt2,Alt4},
which can be applied to the subsonic jet flow problems.
In these works, the free boundary and the solution are
obtained together by solving a minimum problem with free boundary.
Furthermore, the solution on the free boundary must take the extreme value
so that the variational approach works.
Recently, L. L. Du et al \cite{Du1,Du2} used this variational approach
to study impinging subsonic jets and collision of two subsonic flows.
There are also many works on irrotational and rotational subsonic flows past profiles or in nozzles, which are formulated
as fixed boundary problems,
and we refer to \cite{BDX2,Bers1,CX2,DO,DX,DXY,FG2,LY,Weng,XX,YY} and the references therein.
Continuous subsonic-sonic flows in convergent nozzles was studied in \cite{NW1,NW2,W,WX1,WX2},
where the flows are governed by free boundary problems of a degenerate elliptic equation,
and the sonic curve is a free boundary
where the flow velocity is along the normal direction.

It was shown in \cite{Alt1,Alt2,Alt4} that for a nozzle satisfying some assumptions,
if the mass flux of the flow is prescribed, then there exists a unique subsonic jet flow
which is infinitely long and whose pressure on the free boundary is a constant.
Here, the free boundary is the particle path connected with the wall of the nozzle.
Furthermore, the pressure of the subsonic jet flow on the free boundary is a constant
which is determined by the flow to be found.
That is to say, the surrounding pressure cannot be given in advance in these problems.
In many physical problems, the surrounding pressure should be a known constant.
So a natural question is how to formulate a subsonic jet flow whose pressure on the free boundary
coincides with the given surrounding pressure.
In this paper, we study subsonic jet flows for a given surrounding pressure
from a finitely long convergent nozzle with straight solid wall.
At the inlet, the incoming mass flux is a given constant
and the flow angle is prescribed,
which is different from \cite{Alt1,Alt2,Alt4} where the boundary condition at the inlet is
to prescribe the stream function.
If the pressure of the flow in the nozzle is greater than the surrounding pressure,
it is expected that there is an accelerating subsonic jet flow
whose pressure on the free boundary coincides with the given surrounding pressure.
In general, one part of the free boundary should be the particle path
connected with the wall of the nozzle as in \cite{Alt1,Alt2,Alt4}.
Because the surrounding pressure is given in advance,
the subsonic jet flow may be located in a bounded domain.
The outlet of the subsonic jet flow is another part of the free boundary,
and one should prescribe another boundary condition on the outlet
except for the coincidence between the pressure of the flow and the given surrounding pressure.
As in \cite{WX2}, this boundary condition
is prescribed that the flow velocity at the outlet is along the normal direction.

Two-dimensional steady compressible fluids satisfy the Euler system:
\begin{align*}
\pd{}{x}(\rho u)+\pd{}{y}(\rho v)=0,\quad
\pd{}{x}(P+\rho u^2)+\pd{}{y}(\rho uv)=0,\quad\pd{}{x}(\rho uv)+\pd{}{y}(P+\rho v^2)=0,
\end{align*}
where $(u,v)$, $P$ and $\rho$ represent the velocity, pressure and
density of the flow, respectively.
The flow is assumed to be isentropic so that $P=P(\rho)$ is a smooth function. In
particular, for a polytropic gas with adiabatic exponent
$\gamma>1$, $P(\rho)=\rho^\gamma/\gamma$ is the normalized pressure. Assume further that the flow is
irrotational.
Then the density $\rho$ is expressed in terms of the speed
$q$ according to the Bernoulli law (\cite{Bers})
\begin{align}
\label{Bernoulli}
\rho(q^2)=\Big(1-\frac{\gamma-1}{2}q^2\Big)^{1/(\gamma-1)},\quad
0<q<\sqrt{2/{(\gamma-1)}}.
\end{align}
The sound speed $c$ is defined as $c^2=P'(\rho)$.
At the sonic state, it is $c_*=\sqrt{2/{(\gamma+1)}}$,
which is critical in the sense that
the flow is subsonic ($q<c$) when $q<c_*$, sonic ($q=c$) when $q=c_*$,
and supersonic ($q>c$) when $q>c_*$.
It is well-known that the above Euler system can be transformed into the full potential equation (\cite{Bers})
\begin{align}
\label{fullpotential}
\mbox{div}(\rho(|\nabla\varphi|^2)\nabla\varphi)=0,
\end{align}
where $\varphi$ is a velocity potential with $\nabla\varphi=(u,v)$,
and $\rho$ is given by \eqref{Bernoulli}.

Assume that the nozzle is located
symmetrically with respect to the $x$-axis
with the vertex being $(0,0)$ and the angle at the vertex being $2{\vartheta}\in(0,\pi)$.
In this paper we only
consider the upper part of the subsonic jet flow due to the symmetry.
The inlet of the nozzle $\Gamma_{\rm in}$ is the arc centered at $(0,0)$ with radius $R_0>0$.
The upper wall of the nozzle $\Gamma_{\rm w}$ ends at $(-R\cos{\vartheta},R\sin{\vartheta})$ for $R\in(0,R_0)$.
For given constants $P_e>0$ and $m>0$,
we seek a subsonic jet flow in $\Omega$,
which is bounded by $\Gamma_{\rm in}$, $\Gamma_{\rm w}$, the $x$-axis, 
the particle path connected with the upper wall $\Gamma_{\rm ws}$ and
the outlet of the flow $\Gamma_{\rm out}$,
such that the incoming mass flux at $\Gamma_{\rm in}$ is $m$,
the flow velocity at $\Gamma_{\rm in}$ and $\Gamma_{\rm out}$ is along the normal direction,
the flow satisfies the slip condition at $\Gamma_{\rm w}\cup\Gamma_{\rm ws}$ and the $x$-axis,
and the pressure of the flow at $\Gamma_{\rm ws}\cup\Gamma_{\rm out}$ is $P_e$,
where $\Gamma_{\rm ws}$ and $\Gamma_{\rm out}$ are free.
Such a subsonic jet flow problem is formulated as
the following free boundary problem
\begin{align}
\label{phy1}
&\mbox{div}(\rho(|\nabla\varphi|^2)\nabla\varphi)=0,\quad&&(x,y)\in\Omega,
\\
\label{phy3-0}
&-\int_{\Gamma_{\rm in}}\rho(|\nabla\varphi(x,y)|^2)\nabla\varphi(x,y)\cdot\nu(x,y) dl=m,
\\
\label{phy3}
&\varphi(x,y)=0,\quad&&(x,y)\in\Gamma_{\rm in},
\\
\label{phy2}
&\nabla\varphi(x,y)\cdot\nu(x,y)=0,\quad&&(x,y)\in\Gamma_{\rm w}\cup\Gamma_{\rm ws}\cup(\mbox{$x$-axis}\cap\partial\Omega),
\\
\label{phy4}
&P(\rho(|\nabla\varphi(x,y)|^2))=P_e,\quad&&(x,y)\in\Gamma_{\rm ws}\cup\Gamma_{\rm out},
\\
\label{phy5}
&\varphi(x,y)=\xi,\quad&&(x,y)\in\Gamma_{\rm out},
\end{align}
where $\nu$ is the unit outer normal to $\partial\Omega$, $\xi$ is a free constant,
and $\Gamma_{\rm ws}\cup\Gamma_{\rm out}$ is a free boundary.
As to the Dirichlet boundary conditions \eqref{phy3} and \eqref{phy5},
it means that the flow velocity at $\Gamma_{\rm in}$ and $\Gamma_{\rm out}$ is along the normal direction.
Hence both $\Gamma_{\rm in}$ and $\Gamma_{\rm out}$ are level sets of the velocity potential,
and the value at $\Gamma_{\rm in}$ is normalized to be zero without loss of generality.

For the subsonic jet flow problem \eqref{phy1}--\eqref{phy5},
its pressure coincides with the given surrounding pressure at the free boundary.
In general, this free boundary contains two parts: one is the particle path connected with the wall
and the other is a level set of the velocity potential.
Both the velocity potential and the stream function
do not take the extreme value on the whole free boundary.
Hence the variational approach by H. W. Alt, L. A. Caffarelli and A. Friedman
cannot be applied to this free boundary problem.
As far as we know, there are no studies on such problems
although elliptic free boundary problems have been studied extensively
(see, e.g., \cite{CS,Friedman,Lin} and the references therein).
Indeed, it is hard to solve the problem \eqref{phy1}--\eqref{phy5} in the physical plane
since the characteristics of the two parts of the free boundary are completely different.
It is noted that $\Gamma_{\rm ws}$ and $\Gamma_{\rm out}$ are two segments
in the potential-stream coordinates.
So, we study the problem \eqref{phy1}--\eqref{phy5} in the potential plane
since the shape but not the precise location of the free boundary is known in advance
in the potential plane.
However, for the subsonic jet flow problem in the potential plane,
the flow satisfies a nonlinear Robin boundary condition at the inlet
which causes a crucial difficulty for its well-posedness as in \cite{WX2}.
We need to choose a suitable space of solutions to ensure its well-posedness.
Besides this difficulty, there are three new ones completely different from \cite{WX2}.
First of all,
the subsonic jet flow problem in this paper is not a perturbed problem and there are no background solutions,
while \cite{WX2} concerns the structural stability of a symmetric flow.
The other two new difficulties are that
the free boundary contains two different parts,
and mixed Dirichlet-Neumann boundary conditions are prescribed on a segment.
For these difficulties, we need some new estimates and techniques completely different from \cite{WX2},
such as the optimal H\"older estimates for subsonic jet flows,
the continuous dependence of subsonic jet flows with respect to the free boundary,
the precise properties of the free boundary.

To solve the subsonic jet flow problem in the potential plane,
we first study the fixed boundary problem,
for which mixed Dirichlet-Neumann boundary conditions are prescribed on a segment.
We use a duality argument to show the uniqueness of the solution to the fixed boundary problem.
For the dual problem, mixed Dirichlet-Neumann boundary conditions are prescribed
on a segment, and it may be ill-posed due to the boundary condition at the inlet.
We prove the dual problem is well-posed under some restrictions on the
upper bound for the velocity potential difference
and the lower bound for the speed.
As to the existence of solutions to the fixed boundary problem,
we prescribe a Neumann boundary condition instead of the nonlinear Robin one at
the inlet and then use a fixed point argument to get solutions.
To do so, we also need some restrictions on the
upper bound for the velocity potential difference
and the lower bound for the speed.
Since mixed Dirichlet-Neumann boundary conditions are prescribed on a segment,
the regularity of solutions is weak at the joint point of the Dirichlet and Neumann data.
Although there are many studies on mixed boundary value problems for elliptic equations
(see, e.g., \cite{AK,EF,KM,Lieberman1,Lieberman3,Liebermanbook} and the references therein),
the optimal regularity of solutions to this problem is unknown yet.
In this paper, we use a series of elaborate estimates to get the optimal H\"older continuity of solutions.
Summing up, by restricting a suitable upper bound for the velocity potential difference,
we show the well-posedness and the optimal regularity of solutions to the fixed boundary problem
in a suitable space of solutions where there is a lower bound
for the speed.
With the well-posedness of the fixed boundary problem,
together with the continuous dependence of its solutions,
we can solve the subsonic jet flow problem in the potential plane.
We identify a suitable space of flows in terms of
the minimal speed and the maximal velocity potential difference,
by which we get a complete classification of the nozzles
whether there are subsonic jet flows or not.
By studying the precise properties of the subsonic jet flows,
one can transform them into the physical plane.
The main result of this paper is that:
For $P_e=P(\rho(c_e^2))\in(P(\rho(c_*^2)),1/\gamma)$ and ${R_0}{{\vartheta}} {c_l}\rho(c_l^2)<m<{R_0}{{\vartheta}} c_e\rho(c_e^2)$
($c_l\in(0,c_e)$ is a constant),
there exists a positive constant $R_*\in(\hat R,R_0]$
with $\hat R={m}/{({\vartheta}c_e\rho(c_e^2))}$,
such that the problem \eqref{phy1}--\eqref{phy5} admits a unique solution $(\varphi,\Omega)$
in the suitable space if $R\in[\hat R,R_*]\cap[\hat R,R_0)$,
while there is not such a solution if $R\in[0,\hat R)\cup(R_*,R_0)$.
It is noted that the bounds of the incoming mass flux are needed.
For the regularity of the subsonic jet flow, it is proved that $\varphi\in C^{1,\alpha}(\overline\Omega)$
and $\Gamma_{\rm w}\cup\Gamma_{\rm ws}\in C^{1+\alpha}$ for each exponent $\alpha\in(0,1/2)$,
which are almost optimal.
As to the geometry of the free boundary, it is shown that
both $\Gamma_{\rm ws}$ and $\Gamma_{\rm out}$ are strictly convex,
whose tangent lines are located at the same side as the flow.
In this paper, for a class of nozzles, a given surrounding pressure and a given incoming mass flux,
we solve a subsonic jet flow problem in a suitable space,
which is identified for the well-posedness of the problem,
and we get a complete classification of the nozzles
whether there are subsonic flows in this space or not.
However, there may be other subsonic jet flows not in this space,
which will be dealt with in our forthcoming study.
In particular, it is shown that if the angle of the nozzle is small, 
there is an infinite long subsonic or sonic jet flow whose free boundary is particle path
connected the wall and the infinity.

The paper is arranged as follows. In $\S\, 2$, we state the main results of the paper
and formulate the subsonic jet flow problem into a free boundary problem in the potential plane.
It is proved in $\S\,3$ that the fixed boundary problem is well-posed.
Subsequently, the free boundary problem is solved in $\S\, 5$.

\section{Main results and formulation in the potential plane}

In this section, we first state the main results of the paper (well-posedness, nonexistence and properties of solutions).
Then we formulate the subsonic jet flow problem
in the potential plane and introduce the spaces of solutions.

\subsection{Main results}

\begin{definition}
For $P_e>0$, $m>0$ and $0<R<R_0$,
$(\varphi,\Omega)$ is said to be a solution to the free boundary problem \eqref{phy1}--\eqref{phy5},
if $\Gamma_{\rm ws},\Gamma_{\rm out}\in C^1$ such that $\Gamma_{\rm w}\cap\Gamma_{\rm ws}=(-R\cos{\vartheta},R\sin{\vartheta})$
and $\Gamma_{\rm out}$ connects the $x$-axis and $\Gamma_{\rm ws}$,
and $\varphi\in C^2(\overline\Omega\setminus\{\Gamma_{\rm w}\cap\Gamma_{\rm ws},\Gamma_{\rm ws}\cap\Gamma_{\rm out}\})
\cap H^2(\Omega)\cap C^1(\overline\Omega)$
with $0<\inf_{\Omega}|\nabla\varphi|\le\sup_{\Omega}|\nabla\varphi|<c_*$
such that \eqref{phy1}--\eqref{phy5} hold, where
$\Omega$ is the domain bounded by
$\Gamma_{\rm in}$, $\Gamma_{\rm w}$, the $x$-axis, $\Gamma_{\rm ws}$ and $\Gamma_{\rm out}$.
\end{definition}

The main results of the paper are the following theorems.

\begin{theorem}
\label{existencephysical}
For $P(\rho(c_*^2))<P_e<1/\gamma$ and ${R_0}{{\vartheta}} {c_l}\rho(c_l^2)<m<{R_0}{{\vartheta}} c_e\rho(c_e^2)$,
there exists a constant $R_*\in(\hat R,R_0]$
depending only on $R_0$, ${\vartheta}$, $P_e$, $m$ and $\gamma$,
such that the problem \eqref{phy1}--\eqref{phy5} admits a unique solution $(\varphi_{[R]},\Omega_{[R]})$
with $\inf_{\Omega_{[R]}}|\nabla\varphi_{[R]}|\ge {c_l}$ and $\sup_{\Omega_{[R]}}\varphi_{[R]}\le R_0{c_l}$
if $R\in[\hat R,R_*]\cap[\hat R,R_0)$,
while there is not such a solution if $R\in[0,\hat R)\cup(R_*,R_0)$,
where $c_e\in(0,c_*)$ such that $P(\rho(c_e^2))=P_e$,
$\hat R={m}/{({\vartheta}c_e\rho(c_e^2))}$,
and ${c_l}\in(0,c_e)$ is the root to
\begin{align*}
\rho(c_l^2)\int_{c_l}^{c_e}\frac{\rho(q^2)+2q^2\rho'(q^2)}{q\rho^2(q^2)}dq=1.
\end{align*}
\end{theorem}

\begin{remark}
The bounds of the incoming mass flux and restrictions of the solution in Theorem \ref{existencephysical} 
are needed in this paper (see the argument in $\S$ 2.3 and $\S$ 2.4 below).
Furthermore, Theorem \ref{existencephysical} gives a complete classification of the nozzles
whether there are subsonic flows in this space or not.
\end{remark}

\begin{theorem}
\label{existencephysical2}
Assume that $P(\rho(c_*^2))<P_e<1/\gamma$ and ${R_0}{{\vartheta}} {c_l}\rho(c_l^2)<m<{R_0}{{\vartheta}} c_e\rho(c_e^2)$,
$(\varphi_{[R]},\Omega_{[R]})$
with $\inf_{\Omega_{[R]}}|\nabla\varphi_{[R]}|\ge {c_l}$ and $\sup_{\Omega_{[R]}}\varphi_{[R]}\le R_0{c_l}$
is the solution to the problem \eqref{phy1}--\eqref{phy5} for $R\in[\hat R,R_*]\cap[\hat R,R_0)$.

{\rm(i)}
$\Gamma_{\rm ws}$ and $\Gamma_{\rm out}$
can be regarded as the graphs of $y=W_{[R]}(x)\,(-R\cos{\vartheta}\le x\le x_{[R]})$
and $x=J_{[R]}(y)\,(0\le y\le y_{[R]})$, respectively,
where $W_{[R]}\in C^\infty(-R\cos{\vartheta},x_{[R]})\cap C^{1,1}((-R\cos{\vartheta},x_{[R]}])
\cap C^{1,\alpha}([-R\cos{\vartheta},x_{[R]}])$ for each exponent $\alpha\in(0,1/2)$,
$J_{[R]}\in C^\infty([0,y_{[R]}))\cap C^{1,1}([0,y_{[R]}])$, and
\begin{gather*}
-\tan{\vartheta}<W_{[R]}'(x)<0,\quad W_{[R]}''(x)>0,\quad x\in(-R\cos{\vartheta},x_{[R]}),
\\
0<J_{[R]}'(y)<\tan{\vartheta},\quad J_{[R]}''(y)>0,\quad y\in(0, y_{[R]}).
\end{gather*}
Moreover, $\Gamma_{\rm ws}=\emptyset$ if and only if $R=\hat R$.

{\rm(ii)}
$\varphi_{[R]}\in C^\infty(\overline\Omega_{[R]}\setminus\{(-R\cos{\vartheta},R\sin{\vartheta}),(x_{[R]},y_{[R]})\})
\cap C^{1,1}(\overline\Omega_{[R]}\setminus\{(-R\cos{\vartheta},R\sin{\vartheta})\})\cap C^{1,\alpha}(\overline\Omega_{[R]})$
for each exponent $\alpha\in(0,1/2)$, and
\begin{gather*}
{c_l}<|\nabla\varphi_{[R]}(x,y)|< c_e,\quad
-\pd{\varphi_{[R]}}{x}(x,y)\tan{\vartheta}<\pd{\varphi_{[R]}}{y}(x,y)<0,
\quad
(x,y)\in\Omega_{[R]}.
\end{gather*}

{\rm(iii)}
If $R_*<R_0$ additionally, then $\sup_{\Omega_{[R_*]}}\varphi_{[R_*]}=R_0{c_l}$.

{\rm(iv)}
For each $R_1\in[\hat R,R_*]\cap[\hat R,R_0)$,
\begin{gather*}
\lim_{\stackrel{R\to R_1}{R\in[\hat R,R_*]{\cap}[\hat R,R_0)}}\Omega_{[R]}
=\Omega_{[R_1]},
\\
\lim_{\stackrel{R\to R_1}{R\in[\hat R,R_*]{\cap}[\hat R,R_0)}}\nabla\varphi_{[R]}(x,y)
=\nabla\varphi_{[R_1]}(x,y)\mbox{ uniformly for } (x,y)\in\Omega_{[R_1]}.
\end{gather*}

{\rm(v)}
For  $R_1,R_2\in[\hat R,R_*]\cap[\hat R,R_0)$ with
$R_1<R_2$, $\sup_{\Omega_{[R_1]}}\varphi_{[R_1]}<\sup_{\Omega_{[R_2]}}\varphi_{[R_2]}$.
\end{theorem}

\begin{remark}
$\alpha\in(0,1/2)$ in Theorem \ref{existencephysical2} is almost optimal.
\end{remark}

\subsection{Formulation in the potential plane}

Define a velocity potential $\varphi$ and a stream
function $\psi$, respectively, by
$$
\pd\varphi x=u=q\cos\theta,\quad\pd\varphi y=v=q\sin\theta,\quad
\pd\psi x=-\rho v=-\rho q\sin\theta,\quad\pd\psi y=\rho u=\rho q\cos\theta,
$$
where $\theta$ is the flow angle.
The full potential equation \eqref{fullpotential}
can be reduced to the following Chaplygin equations (\cite{Bers}):
\begin{align}
\label{Chaplygin}
\pd\theta\psi+\frac{\rho(q^2)+2q^2\rho'(q^2)}{q\rho^2(q^2)}
\pd{q}\varphi=0,
\quad\frac{1}{q}\pd{q}\psi-\frac{1}{\rho(q^2)}\pd\theta\varphi=0
\end{align}
in the potential-stream coordinates $(\varphi,\psi)$.
And the coordinate transformations between the two coordinate systems
are valid at least in the absence of stagnation points.
Eliminating $\theta$ from \eqref{Chaplygin} yields the following second-order quasilinear equation
\begin{align}
\label{strong}
\pd{^2A(q)}{\varphi^2}+\pd{^2B(q)}{\psi^2}=0,
\end{align}
where
\begin{align*}
A(q)=\int_{c_*}^q\frac{\rho(s^2)+2s^2\rho'(s^2)}{s\rho^2(s^2)}ds,\quad
B(q)=\int_{c_*}^q\frac{\rho(s^2)}{s}ds, \qquad 0<q<\sqrt{2/{(\gamma-1)}}.
\end{align*}
Here, $B(\cdot)$ is strictly increasing in $(0,\sqrt{2/{(\gamma-1)}})$,
while $A(\cdot)$ is strictly increasing in $(0,c_*]$ and
strictly decreasing in $[c_*,\sqrt{2/{(\gamma-1)}})$.
We use $A^{-1}(\cdot)$ to denote the inverse function of $A(\cdot)\big|_{(0,c_*)}$ in this paper.

Assume that $(-R_0,0)$ in the physical plane is transformed into the origin in the potential plane
without loss of generality.
Then $\Omega$ is transformed into $(0,\xi)\times(0,m)$.
Rewrite $\Gamma_{\rm in}$ as
$$
x(s)=-R_0\cos\frac{s}{R_0},\quad y(s)=R_0\sin\frac{s}{R_0},\quad s\in[0,R_0{\vartheta}],
$$
where $s$ is the arc length of $\Gamma_{\rm in}$.
Denote the coordinate transformation from $\{0\}\times[0,m]$ to $\Gamma_{\rm in}$
by $S_{\rm in}(\psi)$. Then $S_{\rm in}(0)=0$, $S_{\rm in}(m)=R_0{\vartheta}$, and
$S'_{\rm in}(\psi)={1}/{(q(0,\psi)\rho(q^2(0,\psi)))}$ for $\psi\in[0,m]$.
Hence
$$
\int_{0}^{m}\frac1{q(0,\psi)\rho(q^2(0,\psi))}d\psi={R_0}{\vartheta}.
$$
It follows from the first equation in \eqref{Chaplygin} that
$$
\pd{A(q)}{\varphi}(0,\psi)=-\pd{\theta}{\psi}(0,\psi)
=\frac{1}{R_0}S'_{\rm in}(\psi)=\frac{1}{{R_0}q(0,\psi)\rho(q^2(0,\psi))},\quad
\psi\in(0,m).
$$
Assume that the velocity potential at $\Gamma_{\rm w}\cap\Gamma_{\rm ws}$ is $\zeta$.
It follows from \eqref{phy2} and the second equation in \eqref{Chaplygin} that
$\pd{B(q)}{\psi}(\cdot,0)\Big|_{(0,\xi)}=\pd{B(q)}{\psi}(\cdot,m)\Big|_{(0,\zeta)}=0$.
Furthermore, \eqref{phy4} yields
$q(\cdot,m)\big|_{(\zeta,\xi)}=q(\xi,\cdot)\big|_{(0,m)}=c_e$,
where $c_e\in(0,c_*)$ such that $P(\rho(c_e^2))=P_e$.
Therefore, the subsonic jet flow problem \eqref{phy1}--\eqref{phy5} is formulated in the potential plane as
the following free boundary problem
\begin{align}
\label{prob1}
&\pd{^2A(q)}{\varphi^2}+\pd{^2B(q)}{\psi^2}=0,
\quad&&(\varphi,\psi)\in(0,\xi)\times(0,m),
\\
\label{prob2}
&\pd{A(q)}{\varphi}(0,\psi)=\frac{1}{{R_0}q(0,\psi)\rho(q^2(0,\psi))},
\quad&&\psi\in(0,m),
\\
\label{prob3}
&\pd{B(q)}{\psi}(\varphi,0)=0,\quad&&\varphi\in(0,\xi),
\\
\label{prob4}
&\pd{B(q)}{\psi}(\varphi,m)=0,\quad&&\varphi\in(0,\zeta),
\\
\label{prob5}
&q(\varphi,m)=c_e,\quad&&\varphi\in(\zeta,\xi),
\\
\label{prob6}
&q(\xi,\psi)=c_e,\quad&&\psi\in(0,m),
\\
\label{prob7}
&\int_{0}^{m}
\frac1{q(0,\psi)\rho(q^2(0,\psi))}d\psi={R_0}{\vartheta},
\end{align}
where $0<c_e<c_*$, $m>0$ and $\zeta>0$ are given constants, while $(q,\xi)$ is the solution.

\begin{definition}
\label{definition}
For $0<c_e<c_*$, $m>0$ and $\zeta>0$,
$(q,\xi)$ is said to be a solution to the free boundary problem \eqref{prob1}--\eqref{prob7},
if $q\in C^2([0,\xi]\times[0,m]\setminus\{(\zeta,m),(\xi,m)\})$
with $0<\inf_{(0,\xi)\times(0,m)}q\le\sup_{(0,\xi)\times(0,m)}q<c_*$
and $\xi\ge\zeta$ such that \eqref{prob1}--\eqref{prob7} hold.
\end{definition}

Solutions, subsolutions and supersolutions to the fixed boundary problem \eqref{prob1}--\eqref{prob6}
can be defined similarly.
Although this fixed boundary problem is uniformly elliptic, the well-posedness, regularity and continuous
dependence of solutions are unknown yet.

\subsection{Bounds of the incoming mass flux}

If $\xi=\zeta$, then the problem \eqref{prob1}--\eqref{prob7} is symmetric, and it can be simplified into
\begin{align}
\label{sprob1}
&\frac{d^2}{d\varphi^2}A(\hat q(\varphi))=0,
\quad\varphi\in(0,\zeta),
\\
\label{sprob2}
&\frac{d}{d\varphi}A(\hat q(0))=\frac{1}{{R_0}\hat q(0)\rho(\hat q^2(0))},
\quad
\hat q(\zeta)=c_e,
\\
\label{sprob4}
&\frac{m}{\hat q(0)\rho(\hat q^2(0))}={R_0}{\vartheta},
\end{align}
where $\zeta>0$ is free.
The free boundary problem \eqref{sprob1}--\eqref{sprob4} admits a solution if and only if $0<m<{R_0}{{\vartheta}} c_e\rho(c_e^2)$.

\begin{lemma}
\label{symmetricexistence}
For $0<c_e<c_*$ and $0<m<{R_0}{{\vartheta}} c_e\rho(c_e^2)$,
the free boundary problem \eqref{sprob1}--\eqref{sprob4} admits a unique solution $\hat q$
with $\zeta=\hat\zeta={m (A(c_e)-A({c_m}))}/{\vartheta}$,
where ${c_m}\in(0,c_e)$ is the unique root to ${R_0}{{\vartheta}} {c_m}\rho(c_m^2)=m$.
More precisely, the solution is
\begin{align*}
\hat q(\varphi)=A^{-1}\big(A(c_e)-{{\vartheta}}(\hat\zeta-\varphi)/m\big),
\quad\varphi\in[0,\hat\zeta].
\end{align*}
In physical plane, it is
\begin{gather*}
\hat\varphi(x,y)=\int^{R_0}_{\sqrt{x^2+y^2}}\hat h\Big(\frac{m}{{\vartheta}r}\Big)dr,\quad(x,y)\in\overline{\hat\Omega},
\\
\hat\Omega=\big\{(x,y)\in\mathbb R^2:\hat R<\sqrt{x^2+y^2}<R_0,x<0,0<y<-x\tan{\vartheta}\big\},
\end{gather*}
where $\hat h$ is the inverse function of $q\rho(q^2)$ in $q\in(0,c_*)$, and
$\hat R={m}/{({\vartheta}c_e\rho(c_e^2))}$.
\end{lemma}

It is noted that $\hat q$ solves the fixed boundary problem \eqref{prob1}--\eqref{prob6}
with $\xi=\zeta=\hat\zeta$.
Assume that $q_1$ and $q_2$ are two solutions to
this problem. Set $Q=A(q_1)-A(q_2)$.
Then $Q$ solves
\begin{align*}
&\pd{^2Q}{\varphi^2}+\pd{^2(b(\varphi,\psi)Q)}{\psi^2}=0,
\quad&&(\varphi,\psi)\in(0,\hat\zeta)\times(0,m),
\\
&\pd{Q}{\varphi}(0,\psi)=-\frac{h(\psi)}{{R_0}}Q(0,\psi),
\quad Q(\xi,\psi)=0,\quad&&\psi\in(0,m),
\\
&\pd{Q}{\psi}(\varphi,0)=0,\quad\pd{Q}{\psi}(\varphi,m)=0,\quad&&\varphi\in(0,\hat\zeta),
\end{align*}
where
\begin{gather*}
b(\varphi,\psi)=\int_0^1\frac{B'}{A'}\big(A^{-1}(tA(q_1(\varphi,\psi))+(1-t)A(q_2(\varphi,\psi)))\big)dt,\quad
(\varphi,\psi)\in(0,\hat\zeta)\times(0,m),
\\
h(\psi)=\int_0^1\frac1{A^{-1}(t{A(q_1(0,\psi))}+(1-t){A(q_2(0,\psi))})}dt,\quad
\psi\in(0,m).
\end{gather*}
Its dual problem is
\begin{align}
\label{acc1}
&\pd{^2U}{\varphi^2}
+b(\varphi,\psi)\pd{^2U}{\psi^2}=0,
\quad&&(\varphi,\psi)\in(0,\hat\zeta)\times(0,m),
\\
&\pd{U}{\varphi}(0,\psi)=-\frac{h(\psi)}{{R_0}}U(0,\psi),
\quad U(\hat\zeta,\psi)=0,\quad&&\psi\in(0,m),
\\
\label{acc2}
&\pd{U}{\psi}(\varphi,0)=0,\quad
\pd{U}{\psi}(\varphi,m)=0,\quad&&\varphi\in(0,\hat\zeta).
\end{align}
If $q_1$ and $q_2$ are small perturbations of $\hat q$, then
$h$ is a small perturbation of $1/{{c_m}}$.
Therefore, the eigenvalue problem for \eqref{acc1}--\eqref{acc2} is
a small perturbation of
\begin{align}
\label{acc3}
&\pd{^2U}{\varphi^2}
+b(\varphi,\psi)\pd{^2U}{\psi^2}+\lambda U=0,
\quad&&(\varphi,\psi)\in(0,\hat\zeta)\times(0,m),
\\
&\pd{U}{\varphi}(0,\psi)=-\frac{1}{{R_0}{c_m}}U(0,\psi),
\quad U(\hat\zeta,\psi)=0,\quad&&\psi\in(0,m),
\\
\label{acc4}
&\pd{U}{\psi}(\varphi,0)=0,\quad
\pd{U}{\psi}(\varphi,m)=0,\quad&&\varphi\in(0,\hat\zeta),
\end{align}
where $\lambda$ is a constant.
If $\hat\zeta\ge{R_0}{c_m}$, it is clear that
the problem \eqref{acc3}--\eqref{acc4} admits a nonpositive eigenvalue.
If $\hat\zeta<{R_0}{c_m}$, one can prove that
there is not a nontrivial solution to the problem \eqref{acc3}--\eqref{acc4}
for each $\lambda\le0$ (the proof can be found
in Proposition \ref{comparison}).
So, to prove that $\hat q$ is the unique solution to
the fixed boundary problem \eqref{prob1}--\eqref{prob6}
with $\xi=\zeta=\hat\zeta$ by a duality argument,
it is reasonable to restrict $\hat\zeta<{R_0}{c_m}$,
which is equivalent to
\begin{align}
\label{m}
{R_0}{{\vartheta}} {c_l}\rho(c_l^2)<m<{R_0}{{\vartheta}} c_e\rho(c_e^2),
\end{align}
where $c_l$ is the constant given in Theorem \ref{existencephysical}.

\begin{remark}
If $\hat\zeta\ge{R_0}{c_m}$, it is unknown whether the solution to the fixed boundary problem \eqref{prob1}--\eqref{prob6}
with $\xi=\zeta=\hat\zeta$ is unique or not.
As mentioned above, $\hat\zeta<{R_0}{c_m}$ is a reasonable restriction
when one proves the uniqueness by a duality argument.
\end{remark}

\begin{remark}
${c_l}\in(0,c_e)$ given in Theorem \ref{existencephysical} satisfies
$\rho(c_l^2)(A(c_e)-A({c_l}))=1$.
Note that $\rho(q^2)(A(c_e)-A(q))$ is strictly decreasing in $q\in(0,c_e)$
and $\lim_{q\to0^+}\rho(q^2)(A(c_e)-A(q))=+\infty$. Hence ${c_l}$ is well defined.
\end{remark}

\subsection{Spaces of solutions to the fixed and free boundary problems}

In order to solve the free boundary problem \eqref{prob1}--\eqref{prob7},
we first show the well-posedness of the fixed boundary problem \eqref{prob1}--\eqref{prob6}
and then determine the free boundary by \eqref{prob7}.

The uniqueness of the solution to the fixed boundary problem \eqref{prob1}--\eqref{prob6}
will be proved by a duality argument in the paper.
As shown in $\S$ 2.3, $m$ should satisfy \eqref{m}. For such $m$, the solution $\hat q$ to
the free boundary problem \eqref{sprob1}--\eqref{sprob4} satisfies $\inf_{(0,\hat\zeta)}\hat q>c_l$.
Hence we choose
\begin{align*}
{\mathscr J}=\Big\{q\in C([0,\xi]\times[0,m]):c_l\le\inf_{(0,\xi)\times(0,m)}q
\le\sup_{(0,\xi)\times(0,m)}q<c_*\Big\}
\end{align*}
as a space of solutions to the problem \eqref{prob1}--\eqref{prob6}.
A similar duality argument as in $\S$ 2.3 shows that
a sufficient condition for the uniqueness of the solution to the problem \eqref{prob1}--\eqref{prob6}
in ${\mathscr J}$ is
\begin{align}
\label{ahatcz}
0<\zeta<R_0 c_l,\quad \zeta\le\xi\le R_0 c_l.
\end{align}
The existence of solutions in ${\mathscr J}$ to the problem \eqref{prob1}--\eqref{prob6}
will be proved by a fixed point argument
as follows: For a given ${g}\in C^{1,\alpha}([0,m])$ satisfying $c_l\le\inf_{(0,m)}{g}\le\sup_{(0,m)}{g}<c_*$,
we solve 
\begin{align}
\label{aprob1}
&\pd{^2A(q)}{\varphi^2}+\pd{^2B(q)}{\psi^2}=0,
\quad&&(\varphi,\psi)\in(0,\xi)\times(0,m),
\\
\label{aprob2}
&\pd{A(q)}{\varphi}(0,\psi)=\frac{1}{{R_0}{g}(\psi)\rho({g}^2(\psi))},
\quad q(\xi,\psi)=c_e,\quad&&\psi\in(0,m),
\\
\label{aprob3}
&\pd{B(q)}{\psi}(\varphi,0)=0,\quad&&\varphi\in(0,\xi),
\\
\label{aprob4}
&\pd{B(q)}{\psi}(\varphi,m)=0,\quad&&\varphi\in(0,\zeta),
\\
\label{aprob6}
&q(\varphi,m)=c_e,\quad&&\varphi\in(\zeta,\xi),
\end{align}
and then define a mapping by ${\mathscr T}({g})=q(0,\cdot)\big|_{[0,m]}$.
The problem \eqref{prob1}--\eqref{prob6} admits a solution if ${\mathscr T}$
has a fixed point.
It is noted that
$$
\overline q(\varphi,\psi)=c_e,\quad
\underline q(\varphi,\psi)=A^{-1}\Big(A(c_e)
-\frac{\xi-\varphi}{{R_0}c_l\rho(c_l^2)}\Big),
\quad (\varphi,\psi)\in[0,\xi]\times[0,m]
$$
are super and sub solutions to the problem \eqref{aprob1}--\eqref{aprob6}, respectively.
If $\zeta$ and $\xi$ satisfy \eqref{ahatcz},
the comparison principle yields
$$
c_l\le A^{-1}\Big(A(c_e)-\frac{\xi}{{R_0}c_l\rho(c_l^2)}\Big)\le({\mathscr T}({g}))(\psi)\le c_e,
\quad \psi\in [0,m].
$$
By this estimate and other ones, one can prove that ${\mathscr T}$ admits a fixed point.

Summing up, we will show the well-posedness of the fixed boundary problem \eqref{prob1}--\eqref{prob6}
in ${\mathscr J}$
for $\zeta$ and $\xi$ satisfying \eqref{ahatcz}.
It can be checked that \eqref{ahatcz} holds for the symmetric flow in
Lemma \ref{symmetricexistence} if and only if $m$ satisfies \eqref{m}.
Hence, we will solve the free boundary problem \eqref{prob1}--\eqref{prob7} in the space
\begin{align*}
{\mathscr S}=
\Big\{(q,\xi):\xi\in(0,R_0c_l],
q\in C([0,\xi]\times[0,m]),
c_l\le\inf_{(0,\xi)\times(0,m)}q
\le\sup_{(0,\xi)\times(0,m)}q<c_*\Big\}
\end{align*}
for ${R_0}{{\vartheta}} {c_l}\rho(c_l^2)<m<{R_0}{{\vartheta}} c_e\rho(c_e^2)$
and $0<\zeta<R_0 c_l$.

\section{Well-posedness of the fixed boundary problem}

In this section, we prove the well-posedness of the fixed boundary problem \eqref{prob1}--\eqref{prob6}
in ${\mathscr J}$,
where $0<c_e<c_*$, $m>0$, and $\zeta$ and $\xi$ satisfy \eqref{ahatcz}.

\subsection{Linear elliptic problem with mixed Dirichlet-Neumann boundary conditions}

Assume that $m>0$, $0<\zeta\le\xi$, $b\in C^1((0,\xi)\times(0,m))\cap H^1((0,\xi)\times(0,m))$
satisfying
\begin{align}
\label{ab}
b_1\le b(\varphi,\psi)\le b_2,\quad
(\varphi,\psi)\in(0,\xi)\times(0,m),
\end{align}
where $b_1\le b_2$ are positive constants.
Consider the following problem
\begin{align}
\label{l-equation}
&\pd{^2U}{\varphi^2}
+b(\varphi,\psi)\pd{^2U}{\psi^2}=\omega(\varphi,\psi),
\quad&&(\varphi,\psi)\in(0,\xi)\times(0,m),
\\
\label{l-bc-0}
&\pd{U}{\varphi}(0,\psi)=h(\psi),
\quad U(\xi,\psi)=0,\quad&&\psi\in(0,m),
\\
\label{l-bc-1}
&\pd{U}{\psi}(\varphi,0)=0,\quad&&\varphi\in(0,\xi),
\\
\label{l-bc-2}
&\pd{U}{\psi}(\varphi,m)=0,\quad&&\varphi\in(0,\zeta),
\\
\label{l-bc-y}
&U(\varphi,m)=0,\quad&&\varphi\in(\zeta,\xi),
\end{align}
where $\omega\in C^1([0,\xi]\times[0,m])$ and $h\in L^\infty(0,m)$.

\begin{definition}
A function $U\in C^2((0,\xi)\times(0,m))\cap H^1((0,\xi)\times(0,m))\cap C([0,\xi]\times[0,m])$ is said to be
a subsolution (supersolution, solution) to the problem \eqref{l-equation}--\eqref{l-bc-y}, if
\begin{align*}
\int_0^\xi\int_0^m \Big(\pd{U}{\varphi}\,\pd{\eta}{\varphi}+\pd{U}{\psi}\,\pd{(b\eta)}{\psi}+\omega\eta\Big)d\varphi d\psi
+\int_0^m h(\psi)\eta(0,\psi)d\psi\le (\ge,=) 0
\end{align*}
for any nonnegative $\eta\in H^1((0,\xi)\times(0,m))\cap C([0,\xi]\times[0,m])$
with $\eta=0$ on $(\zeta,\xi)\times\{m\}\cup\{\xi\}\times(0,m)$,
and $U(\cdot,m)\big|_{(\zeta,\xi)}\le (\ge,=)0$ and $U(\xi,\cdot)\big|_{(0,m)}\le (\ge,=)0$ hold.
\end{definition}

Lieberman has proved the well-posedness of general linear mixed boundary value problems for elliptic equations
in weighted H\"older spaces in \cite{Lieberman1,Lieberman3}.
For the problem \eqref{l-equation}--\eqref{l-bc-y},
we can show its optimal H\"older continuity by suitable super and sub solutions in the following proposition.
Moreover, the solution is still H\"older continuous if
$b\in C([0,\xi]\times[0,m])$ is relaxed to
that the oscillation of $b$ near $(\zeta,m)$ is suitably small (see Proposition \ref{linearpp} below).

\begin{proposition}
\label{linearp}
Assume that $m>0$, $0<\zeta\le\xi$,
$b\in C^1((0,\xi)\times(0,m))\cap H^1((0,\xi)\times(0,m))\cap C([0,\xi]\times[0,m])$
satisfying \eqref{ab}, $\omega\in C^{1}([0,\xi]\times[0,m])$
and $h\in L^\infty(0,m)$.

{\rm (i)}
If $\underline U,\overline U\in C^2((0,\xi)\times(0,m))\cap H^1((0,\xi)\times(0,m))\cap C([0,\xi]\times[0,m])$
are sub and super solutions to the problem \eqref{l-equation}--\eqref{l-bc-y}, respectively, then
$\underline U\le\overline U$ in $(0,\xi)\times(0,m)$.

{\rm (ii)}
The problem \eqref{l-equation}--\eqref{l-bc-y} admits a unique solution
$U\in C^2((0,\xi)\times(0,m))\cap H^1((0,\xi)\times(0,m))\cap C^\alpha([0,\xi]\times[0,m])$
for each exponent $\alpha\in(0,1/2)$. Furthermore,
\begin{gather}
\label{nl3-00}
\|U\|_{C^{\alpha}([0,\xi]\times[0,m])}\le M,
\end{gather}
where $M>0$ is a constant depending only on $m$, $\zeta$, $\xi$,
$\alpha$, $b$, $\|\omega\|_{L^\infty((0,\xi)\times(0,m))}$ and $\|h\|_{L^\infty(0,m)}$.
\end{proposition}

\Proof
It is clear that the comparison principle holds for the problem \eqref{l-equation}--\eqref{l-bc-y}.
As to the existence, it is assumed that $\xi>\zeta$. The proof for the case $\xi=\zeta$ is simpler.
For each positive integer $n$, choose $h_n\in C^2([0,m])$ such that
$\|h_n-h\|_{L^2(0,m)}\le 1/n$ and $\|h_n\|_{L^\infty(0,m)}\le \|h\|_{L^\infty(0,m)}$,
and set
$$
f_n(\varphi)=\left\{
\begin{aligned}
&m,&&\mbox{ if }0\le\varphi\le\zeta,
\\
&m-(\delta^4/{n^4}-(\varphi-\zeta-\delta/n)^4)^{1/4}, &&\mbox{ if }\zeta<\varphi<\zeta+\delta/n,
\\
&m-\delta/n, &&\mbox{ if }\zeta+\delta/n\le\varphi\le\xi-\delta/n,
\\
&m-2\delta/n+(\delta^4/{n^4}-(\varphi-\xi+\delta/n)^4)^{1/4}, &&\mbox{ if }\xi-\delta/n<\varphi\le\xi,
\end{aligned}
\right.
$$
where $\delta=1/3\min\{m,\xi-\zeta\}$.
Consider the following problem
\begin{align}
\label{nl-equation}
&\pd{^2U_n}{\varphi^2}
+b(\varphi,\psi)\pd{^2U_n}{\psi^2}=\omega(\varphi,\psi),
\quad&&(\varphi,\psi)\in G_n,
\\
\label{nl-bc-0}
&\pd{U_n}{\varphi}(0,\psi)=h_n(\psi),
\quad&&\psi\in(0,m),
\\
\label{nl-bc-1}
&\pd{U_n}{\psi}(\varphi,0)=0,\quad&&\varphi\in(0,\xi),
\\
\label{nl-bc-2}
&\pd{U_n}{\psi}(\varphi,m)=0,\quad&&\varphi\in(0,\zeta),
\\
\label{nl-bc-3}
&U_n(\varphi,f_n(\varphi))=0,\quad&&\varphi\in(\zeta,\xi),
\\
\label{nl-bc-y}
&U_n(\xi,\psi)=0,\quad&&\psi\in(0,m-2\delta/n),
\end{align}
where $G_n=\{(\varphi,\psi)\in\mathbb R^2:0<\varphi<\xi,0<\psi<f_n(\varphi)\}$.
The classical theory on elliptic
equations (\cite{GT}$\S\,6.7$) yields that the problem \eqref{nl-equation}--\eqref{nl-bc-y}
admits a unique solution $U_n\in C^2(\overline G_n)$ such that
\begin{gather}
\label{nl1}
\|U_n\|_{L^\infty(G_n)}+\|U_n\|_{H^1(G_n)}\le M_1,
\end{gather}
where $M_1>0$ is a constant depending only on $m$, $\zeta$, $\xi$, $b$,
$\|\omega\|_{L^\infty((0,\xi)\times(0,m))}$ and $\|h\|_{L^\infty(0,m)}$.

We now prove that $U_n$ is uniformly H\"older continuous.
For a given exponent $\alpha\in(0,1/2)$,
set
$$
\tilde U_n(\varphi,\psi)=
- r^{\alpha}(\varphi,\psi)\sin{\beta}(\varphi,\psi),
\quad(\varphi,\psi)\in[0,\xi]\times[0,m],
$$
where
$$
r(\varphi,\psi)=
\sqrt{\lambda^2(\varphi-\zeta_n)^2+(\psi-m)^2},\quad
(\varphi,\psi)\in[0,\xi]\times[0,m],
$$
$$
\beta(\varphi,\psi)=\left\{
\begin{aligned}
&\frac1\mu\arctan\frac{\psi-m}{\lambda (\varphi-\zeta_n)}-\frac{\mu-2}{2\mu}\pi,
&&\quad \mbox{ if }(\varphi,\psi)\in(\zeta_n,\xi]\times[0,m],
\\
&-\frac{\mu-1}{2\mu}\pi,
&&\quad \mbox{ if }\varphi=\zeta_n,\,\psi\in[0,m],
\\
&\frac1\mu\arctan\frac{\psi-m}{\lambda (\varphi-\zeta_n)}-\frac{1}{2}\pi,
&&\quad \mbox{ if }(\varphi,\psi)\in[0,\zeta_n)\times[0,m],
\end{aligned}
\right.
$$
$\zeta_n=\zeta+\delta/n$, $\lambda=\sqrt{b(\zeta,m)}$ and $\mu=1+1/{(2\alpha)}$.
Then $\tilde U_n\in C^{2}(\overline G_n)$
and for $(\varphi,\psi)\in \overline G_n$,
\begin{gather*}
\pd{\tilde U_n}{\varphi}=-\alpha\lambda^2(\varphi-\zeta_n)r^{\alpha-2}\sin{\beta}
+\frac{\lambda}{\mu}(\psi-m)r^{\alpha-2}\cos{\beta},
\\
\pd{\tilde U_n}{\psi}=-\alpha (\psi-m)r^{\alpha-2}\sin{\beta}
-\frac{\lambda}{\mu}(\varphi-\zeta_n)r^{\alpha-2}\cos{\beta},
\end{gather*}
\vskip-6mm
\begin{align*}
\pd{^2\tilde U_n}{\varphi^2}=&
-\alpha(\alpha-2)\lambda^4(\varphi-\zeta_n)^2r^{\alpha-4}\sin{\beta}
-\alpha \lambda^2r^{\alpha-2}\sin{\beta}
+\frac{\alpha\lambda^3}{\mu}(\varphi-\zeta_n)(\psi-m)r^{\alpha-4}\cos{\beta}
\\
&\qquad+\frac{(\alpha-2)\lambda^3}{\mu}(\varphi-\zeta_n)(\psi-m)r^{\alpha-4}\cos{\beta}
+\frac{\lambda^2}{\mu^2}(\psi-m)^2r^{\alpha-4}\sin{\beta},
\end{align*}
\vskip-6mm
\begin{align*}
\pd{^2\tilde U_n}{\psi^2}=&
-\alpha(\alpha-2)(\psi-m)^2r^{\alpha-4}\sin{\beta}
-\alpha r^{\alpha-2}\sin{\beta}
-\frac{\alpha\lambda}{\mu}(\varphi-\zeta_n)(\psi-m)r^{\alpha-4}\cos{\beta}
\\
&\qquad-\frac{(\alpha-2)\lambda}{\mu}(\varphi-\zeta_n)(\psi-m)r^{\alpha-4}\cos{\beta}
+\frac{\lambda^2}{\mu^2}(\varphi-\zeta_n)^2r^{\alpha-4}\sin{\beta}.
\end{align*}
Therefore
\begin{align}
\label{az-1}
&\pd{^2\tilde U_n}{\varphi^2}(\varphi,\psi)
+b(\varphi,\psi)\pd{^2\tilde U_n}{\psi^2}(\varphi,\psi)
\nonumber
\\
=&\Big(\alpha(2-\alpha)b\Big(\frac{\lambda^4}{b}(\varphi-\zeta_n)^2+(\psi-m)^2\Big)
+\frac{\lambda^2}{\mu^2}(b(\varphi-\zeta_n)^2+(\psi-m)^2)
-\alpha(\lambda^2+b)r^2\Big)r^{\alpha-4}\sin{\beta}
\nonumber
\\
&\qquad
+\frac{2(\alpha-1)\lambda}{\mu}(\lambda^2-b)(\varphi-\zeta_n)(\psi-m)r^{\alpha-4}\cos{\beta},
\quad(\varphi,\psi)\in G_n.
\end{align}
Since $\alpha\in(0,1/2)$ and $b\in C([0,\xi]\times[0,m])$ satisfying \eqref{ab},
it follows from the choice of $\mu$ and $\lambda$ that
there exist two positive constants $\tau_1$ and $\tau_2$,
which depend only on $\alpha$ and $b$,
such that
\begin{align}
\label{az-2}
&\alpha(2-\alpha)b\Big(\frac{\lambda^4}{b}(\varphi-\zeta_n)^2+(\psi-m)^2\Big)
+\frac{\lambda^2}{\mu^2}(b(\varphi-\zeta_n)^2+(\psi-m)^2)
\nonumber
\\
\ge&(\alpha(\lambda^2+b)+\tau_2)r^2,
\quad
(\varphi,\psi)\in(0,\xi)\times(0,m),\,(\varphi-\zeta)^2+(\psi-m)^2<\tau_1.
\end{align}
The definition of $\beta$ yields
$\sin\beta\le-\sin(1/2-1/\mu)\pi<0$ in $(0,\xi)\times(0,m)$,
which, together with \eqref{az-1} and \eqref{az-2}, leads to
\begin{align*}
&\pd{^2\tilde U_n}{\varphi^2}(\varphi,\psi)
+b(\varphi,\psi)\pd{^2\tilde U_n}{\psi^2}(\varphi,\psi)
\\
\le&-\tau_2r^{\alpha-2}\sin(1/2-1/\mu)\pi
+\frac{2(1-\alpha)\lambda}{\mu}|(\lambda^2-b)(\varphi-\zeta_n)(\psi-m)|r^{\alpha-4}
\\
\le&-\tau_2r^{\alpha-2}\sin(1/2-1/\mu)\pi
+\frac{1-\alpha}{\mu}|\lambda^2-b|r^{\alpha-2},
\quad (\varphi,\psi)\in G_n,\,
(\varphi-\zeta)^2+(\psi-m)^2<\tau_1.
\end{align*}
From $b\in C([0,\xi]\times[0,m])$ and the choice of $\lambda$,
there exists two constants $\tau_3\in(0,\tau_1]$ and $\tau_4>0$, which depend only on $\alpha$ and $b$,
such that
\begin{align}
\label{ass-5}
\pd{^2\tilde U_n}{\varphi^2}(\varphi,\psi)
+b(\varphi,\psi)\pd{^2\tilde U_n}{\psi^2}(\varphi,\psi)\le-\tau_4 r^{\alpha-2},
\quad (\varphi,\psi)\in G_n,\,
(\varphi-\zeta)^2+(\psi-m)^2<\tau_3.
\end{align}
Using the comparison principle, together with
\eqref{ass-5} and \eqref{nl1}, one gets that
\begin{align}
\label{nl2}
|U_n(\varphi,\psi)|\le M_2\tilde U_n(\varphi,\psi),
\quad (\varphi,\psi)\in G_n,\,
(\varphi-\zeta)^2+(\psi-m)^2<\tau,
\end{align}
where $M_2>0$ and $\tau>0$ are suitable constants depending only on
$m$, $\zeta$, $\xi$, $\alpha$, $b$, $\|\omega\|_{L^\infty((0,\xi)\times(0,m))}$ and $\|h\|_{L^\infty(0,m)}$.
From \eqref{nl1} and \eqref{nl2}, one can prove that
the problem \eqref{l-equation}--\eqref{l-bc-y} admits a solution
$U\in H^1((0,\xi)\times(0,m))$ satisfying 
\begin{align}
\label{nl2-01}
|U(\varphi,\psi)|\le M_3((\varphi-\zeta)^2+(\psi-m)^2)^{\alpha/2},
\quad (\varphi,\psi)\in (0,\xi)\times(0,m),
\end{align}
where $M_3>0$ is a constant depending only on
$m$, $\zeta$, $\xi$, $\alpha$, $b$, $\|\omega\|_{L^\infty((0,\xi)\times(0,m))}$ and $\|h\|_{L^\infty(0,m)}$.
From \eqref{nl2-01}, the Schauder theory and the H\"older estimates for elliptic equations
(\cite{GT}$\S\, 6.7 \& 8.10$), one can get that 
$U\in C^2((0,\xi)\times(0,m))\cap C^{\alpha}([0,\xi]\times[0,m])$ satisfying \eqref{nl3-00}.
$\hfill\Box$\vskip 4mm

\begin{remark}
If $b$ is a positive constant, then $\tilde U_n$ with $\alpha=1/2$ in the proof of Proposition \ref{linearp}
solves the homogeneous equation of \eqref{l-equation}.
Hence the H\"older continuity in Proposition \ref{linearp}
is almost optimal.
\end{remark}

\begin{proposition}
\label{linearpp}
Assume that $m>0$, $0<\zeta\le\xi$,
$b\in C^1((0,\xi)\times(0,m))\cap H^1((0,\xi)\times(0,m))$
satisfying \eqref{ab}. There exist an exponent $\alpha\in(0,1/2)$
and a constant $\sigma>0$, depending only on ${b_1}$ and ${b_2}$,
such that if the oscillation of $b$ near $(\zeta,m)$ is not greater than $\sigma$,
then for $\omega\in C^{1}([0,\xi]\times[0,m])$ and $h\in L^\infty(0,m)$,
the problem \eqref{l-equation}--\eqref{l-bc-y} admits a unique solution
$U\in C^2((0,\xi)\times(0,m))\cap H^1((0,\xi)\times(0,m))\cap C^{\alpha}([0,\xi]\times[0,m])$.
\end{proposition}

\Proof
The proof is similar to Proposition \ref{linearp}
and one needs only to construct a suitable supersolution to the problem \eqref{nl-equation}--\eqref{nl-bc-y}.
Set
$$
\tilde U_n(\varphi,\psi)=
- r^{\alpha}(\varphi,\psi)\sin{\beta}(\varphi,\psi),
\quad(\varphi,\psi)\in[0,\xi]\times[0,m],
$$
where $\alpha\in(0,1)$ is a constant to be determined,
$$
r(\varphi,\psi)=
\sqrt{\lambda^2(\varphi-\zeta_n)^2+(\psi-m)^2},\quad
(\varphi,\psi)\in[0,\xi]\times[0,m],
$$
$$
\beta(\varphi,\psi)=\left\{
\begin{aligned}
&\frac13\arctan\frac{\psi-m}{\lambda (\varphi-\zeta_n)}-\frac{1}{6}\pi,
&&\quad \mbox{ if }(\varphi,\psi)\in(\zeta_n,\xi]\times[0,m],
\\
&-\frac{1}{3}\pi,
&&\quad \mbox{ if }\varphi=\zeta_n,\,\psi\in[0,m],
\\
&\frac13\arctan\frac{\psi-m}{\lambda (\varphi-\zeta_n)}-\frac{1}{2}\pi,
&&\quad \mbox{ if }(\varphi,\psi)\in[0,\zeta_n)\times[0,m],
\end{aligned}
\right.
$$
$\zeta_n=\zeta+\delta/n$ and $\lambda=\sqrt{b(\zeta,m)}$.
Direct calculations give that for $(\varphi,\psi)\in \overline G_n$,
\begin{align*}
&\pd{^2\tilde U_n}{\varphi^2}(\varphi,\psi)
+b(\varphi,\psi)\pd{^2\tilde U_n}{\psi^2}(\varphi,\psi)
\\
=&\Big(\alpha(2-\alpha)({\lambda^4}(\varphi-\zeta_n)^2+b(\psi-m)^2)
+\frac{\lambda^2}{9}(b(\varphi-\zeta_n)^2+(\psi-m)^2)
\\
&\qquad-\alpha(\lambda^2+b)r^2\Big)r^{\alpha-4}\sin{\beta}
+\frac{2(\alpha-1)\lambda}{3}(\lambda^2-b)(\varphi-\zeta_n)(\psi-m)r^{\alpha-4}\cos{\beta}
\\
\le&-\Big(\alpha(2-\alpha)(b_1^2(\varphi-\zeta_n)^2+b_1(\psi-m)^2)
+\frac{b_1}{9}(b_1(\varphi-\zeta_n)^2+(\psi-m)^2)
\\
&\qquad-2\alpha b_2(b_2(\varphi-\zeta_n)^2+(\psi-m)^2)\Big)r^{\alpha-4}\sin\frac\pi6
\\
&\qquad
+\frac{1}{3}|b(\zeta,m)-b(\varphi,\psi)|(b_2(\varphi-\zeta_n)^2+(\psi-m)^2)r^{\alpha-4}
\\
\le&-\Big(\frac{\alpha(2-\alpha)b_1^2}{2b_2}
+\frac{b_1^2}{18b_2}-\alpha b_2
-\frac{|b(\zeta,m)-b(\varphi,\psi)|}{3}\Big)(b_2(\varphi-\zeta_n)^2+(\psi-m)^2)r^{\alpha-4}.
\end{align*}
Choose
$\alpha={b_1^2}/{(36b_2^2)}$ and $\sigma={3\alpha(2-\alpha)b_1^2}/{(2b_2)}$.
If there exists a positive constant $\tau$ such that
$$
|b(\zeta,m)-b(\varphi,\psi)|\le\sigma,\quad
(\varphi,\psi)\in (0,\xi)\times(0,m),\,
(\varphi-\zeta)^2+(\psi-m)^2<\tau,
$$
then
\begin{align*}
\pd{^2\tilde U_n}{\varphi^2}(\varphi,\psi)
+b(\varphi,\psi)\pd{^2\tilde U_n}{\psi^2}(\varphi,\psi)
\le-\frac{b_1^2}{36b_2}r^{\alpha-2},\quad
(\varphi,\psi)\in G_n,\,
(\varphi-\zeta)^2+(\psi-m)^2<\tau.
\end{align*}
Subsequently, one can prove the proposition similarly to Proposition \ref{linearp}.
$\hfill\Box$\vskip 4mm

\subsection{Comparison principle}

\begin{proposition}[Comparison principle]
\label{comparison}
Assume that $0<c_e<c_*$, $m>0$, and $\zeta$ and $\xi$ satisfy \eqref{ahatcz}.
Let $q_1,q_2\in C^2([0,\xi]\times[0,m]\setminus\{(\zeta,m),(\xi,m)\})
\cap H^1((0,\xi)\times(0,m))\cap {\mathscr J}$
be sub and super solutions to the problem \eqref{prob1}--\eqref{prob6}, respectively.
Then $q_1\le q_2$ in $(0,\xi)\times(0,m)$.
\end{proposition}

\Proof
The proof is based on a duality argument. Set
\begin{gather*}
b(\varphi,\psi)=\int_0^1\frac{B'}{A'}\big(A^{-1}(tA(q_1(\varphi,\psi))+(1-t)A(q_2(\varphi,\psi)))\big)dt,\quad
(\varphi,\psi)\in(0,\xi)\times(0,m),
\\
h(\psi)=\int_0^1\frac1{A^{-1}(t{A(q_1(0,\psi))}+(1-t){A(q_2(0,\psi))})}dt,\quad
\psi\in(0,m).
\end{gather*}
Then, $b\in C^2([0,\xi]\times[0,m]\setminus\{(\zeta,m),(\xi,m)\})\cap H^1((0,\xi)\times(0,m))
\cap C([0,\xi]\times[0,m])$ satisfying $b_1\le b\le b_2$ in $(0,\xi)\times(0,m)$ with some positive constants $b_1\le b_2$,
and $h\in C^2([0,m])$ satisfying $1/{c_*}\le h\le 1/{c_l}$ in $(0,m)$.

For each nonpositive $\omega\in C^1([0,\xi]\times[0,m])$, consider the problem
\begin{align}
\label{cth-equation}
&\pd{^2U}{\varphi^2}
+b(\varphi,\psi)\pd{^2U}{\psi^2}=\omega(\varphi,\psi),
\quad&&(\varphi,\psi)\in(0,\xi)\times(0,m),
\\
\label{cth-bc-0}
&\pd{U}{\varphi}(0,\psi)=-\frac{h(\psi)}{{R_0}}U(0,\psi),
\quad U(\xi,\psi)=0,\quad&&\psi\in(0,m),
\\
\label{cth-bc-1}
&\pd{U}{\psi}(\varphi,0)=0,\quad&&\varphi\in(0,\xi),
\\
\label{cth-bc-2}
&\pd{U}{\psi}(\varphi,m)=0,\quad&&\varphi\in(0,\zeta),
\\
\label{cth-bc-y}
&U(\varphi,m)=0,\quad&&\varphi\in(\zeta,\xi).
\end{align}
We prove the well-posedness of the problem \eqref{cth-equation}--\eqref{cth-bc-y}
by the contraction mapping principle.
Set
${\mathscr C}=\big\{{\mathscr U}\in C([0,m]):{\mathscr U}\ge0\mbox{ in }(0,m)\big\}$.
For each ${\mathscr U}\in{\mathscr C}$, it follows from Proposition \ref{linearp} that the problem
\begin{align}
\label{lcth-equation}
&\pd{^2U}{\varphi^2}
+b(\varphi,\psi)\pd{^2U}{\psi^2}=\omega(\varphi,\psi),
\quad&&(\varphi,\psi)\in(0,\xi)\times(0,m),
\\
\label{lcth-bc-0}
&\pd{U}{\varphi}(0,\psi)=-\frac{h(\psi)}{{R_0}}{\mathscr U}(\psi),
\quad U(\xi,\psi)=0,\quad&&\psi\in(0,m),
\\
\label{lcth-bc-1}
&\pd{U}{\psi}(\varphi,0)=0,\quad&&\varphi\in(0,\xi),
\\
\label{lcth-bc-2}
&\pd{U}{\psi}(\varphi,m)=0,\quad&&\varphi\in(0,\zeta),
\\
\label{lcth-bc-y}
&U(\varphi,m)=0,\quad&&\varphi\in(\zeta,\xi)
\end{align}
admits a unique solution $0\le U\in C^2((0,\xi)\times(0,m))\cap H^1((0,\xi)\times(0,m))\cap C([0,\xi]\times[0,m])$.
Therefore, we can define a mapping $J$ from ${\mathscr C}$ to itself
by $J({\mathscr U})=U(0,\cdot)\big|_{[0,m]}$.
For ${\mathscr U}_1,{\mathscr U}_2\in{\mathscr C}$, one has
$J({\mathscr U}_1)-J({\mathscr U}_2)=\tilde U(0,\cdot)\big|_{[0,m]}$
with $\tilde U\in C^2((0,\xi)\times(0,m))\cap H^1((0,\xi)\times(0,m))\cap C([0,\xi]\times[0,m])$ solving
\begin{align}
\label{clcth-equation}
&\pd{^2\tilde U}{\varphi^2}
+b(\varphi,\psi)\pd{^2\tilde U}{\psi^2}=0,
\quad&&(\varphi,\psi)\in(0,\xi)\times(0,m),
\\
\label{clcth-bc-0}
&\pd{\tilde U}{\varphi}(0,\psi)=-\frac{h(\psi)}{{R_0}}({\mathscr U}_1(\psi)-{\mathscr U}_2(\psi)),
\quad \tilde U(\xi,\psi)=0,\quad&&\psi\in(0,m),
\\
\label{clcth-bc-1}
&\pd{\tilde U}{\psi}(\varphi,0)=0,\quad&&\varphi\in(0,\xi),
\\
\label{clcth-bc-2}
&\pd{\tilde U}{\psi}(\varphi,m)=0,\quad&&\varphi\in(0,\zeta),
\\
\label{clcth-bc-y}
&\tilde U(\varphi,m)=0,\quad&&\varphi\in(\zeta,\xi).
\end{align}
It follows from Proposition \ref{linearp} that the problem
\begin{align}
\label{vclcth-equation}
&\pd{^2V}{\varphi^2}
+b(\varphi,\psi)\pd{^2V}{\psi^2}=0,
\quad&&(\varphi,\psi)\in(0,\xi)\times(0,m),
\\
\label{vclcth-bc-0}
&\pd{V}{\varphi}(0,\psi)=-\frac{1}{{R_0}c_l},
\quad V(\xi,\psi)=0,\quad&&\psi\in(0,m),
\\
\label{vclcth-bc-1}
&\pd{V}{\psi}(\varphi,0)=0,\quad&&\varphi\in(0,\xi),
\\
\label{vclcth-bc-2}
&\pd{V}{\psi}(\varphi,m)=0,\quad&&\varphi\in(0,\zeta),
\\
\label{vclcth-bc-y}
&V(\varphi,m)=0,\quad&&\varphi\in(\zeta,\xi)
\end{align}
admits a unique solution $0\le V\in C^2((0,\xi)\times(0,m))\cap H^1((0,\xi)\times(0,m))\cap C([0,\xi]\times[0,m])$.
It is noted that
$$
\tilde U{_\pm}=\pm\|{\mathscr U}_1-{\mathscr U}_2\|_{L^\infty(0,m)}V(\varphi,\psi),
\quad(\varphi,\psi)\in[0,\xi]\times[0,m]
$$
are super and sub solutions to the problem \eqref{clcth-equation}--\eqref{clcth-bc-y}, respectively.
Proposition \ref{linearp} shows that $J$ is a contraction mapping if $V$ satisfies
\begin{align}
\label{v}
\|V\|_{L^\infty((0,\xi)\times(0,m))}<1.
\end{align}
Set
$$
\overline V(\varphi,\psi)=\frac{1}{{R_0}c_l}(\xi-\varphi), \quad(\varphi,\psi)\in[0,\xi]\times[0,m].
$$
Then $\overline V$ is a supersolution to the problem \eqref{vclcth-equation}--\eqref{vclcth-bc-y},
and Proposition \ref{linearp} leads to
\begin{align}
\label{v1}
0\le V(\varphi,\psi)\le\overline V(\varphi,\psi)=\frac{1}{{R_0}c_l}(\xi-\varphi), \quad(\varphi,\psi)\in[0,\xi]\times[0,m].
\end{align}
If $\xi<{R_0}c_l$, then \eqref{v1} yields \eqref{v}.
Turn to the other case that $0<\zeta<\xi=R_0c_l$.
It is noted that both $V$ and $\overline V$ solve \eqref{vclcth-equation}.
The Hopf Lemma yields that
$V(0,\psi)<\overline V(0,\psi)=1$ for $\psi\in[0,m]$,
which and \eqref{v1} imply \eqref{v}.
Summing up, if $\zeta$ and $\xi$ satisfy \eqref{ahatcz},
then $J$ is a contraction mapping.
Therefore, $J$ admits a unique fixed point,
and there exists a unique solution $U\in C^2((0,\xi)\times(0,m))\cap H^1((0,\xi)\times(0,m))\cap C([0,\xi]\times[0,m])$
to the problem \eqref{cth-equation}--\eqref{cth-bc-y}.
Furthermore, it follows from the classical theory on elliptic
equations (\cite{GT}Theorems 6.24 and 6.30) that $U\in C^2([0,\xi]\times[0,m]\setminus\{(\zeta,m),(\xi,m)\})$.

For a positive integer $n$, let $0\le\eta_n\in C^2(\mathbb R^2)$ such that
$$
\eta_n(\varphi,\psi)=\left\{
\begin{aligned}
&0, &&\quad \mbox{ if }(\varphi-\zeta)^2+(\psi-m)^2<\frac1{n^2},
\\
&1, &&\quad \mbox{ if }(\varphi-\zeta)^2+(\psi-m)^2>\frac4{n^2},
\end{aligned}
\right.
$$
\vskip-2mm
$$
|\nabla\eta_n(\varphi,\psi)|\le 4n,\quad
\Big|\pd{^2\eta_n}{\varphi^2}(\varphi,\psi)\Big|+
\Big|\pd{^2\eta_n}{\psi^2}(\varphi,\psi)\Big|\le 8{n^2},
\quad
(\varphi,\psi)\in\mathbb R^2.
$$
It follows from the definition of sub and super solutions that
\begin{align*}
&\int_{0}^{\xi}\int_0^{m}
(A(q_1)-A(q_2))\eta_n\Big(\pd{^2U}{\varphi^2}
+b\pd{^2U}{\psi^2}\Big)d\varphi d\psi
\\
&\qquad+\int_{0}^{\xi}\int_0^{m}
(A(q_1)-A(q_2))U\Big(\pd{^2\eta_n}{\varphi^2}
+b\pd{^2\eta_n}{\psi^2}\Big)d\varphi d\psi
\\
&\qquad+2\int_{0}^{\xi}\int_0^{m}
(A(q_1)-A(q_2))\Big(\pd{\eta_n}\varphi\,\pd U\varphi+b\pd{\eta_n}\psi\,\pd U{\psi}\Big)d\varphi d\psi
\\
\ge&-\int_0^{m}({{A}(q_1)}(0,\psi)-{{A}(q_2)}(0,\psi))
\pd{\eta_n}{\varphi}(0,\psi)U(0,\psi)d\psi
\\
&\qquad-\int_0^{\xi}({{B}(q_1)}(\varphi,0)-{{B}(q_2)}(\varphi,0))
\pd{\eta_n}{\psi}(\varphi,0)U(\varphi,0)d\varphi
\\
&\qquad+\int_0^{\zeta}({{B}(q_1)}(\varphi,m)-{{B}(q_2)}(\varphi,m))
\pd{\eta_n}{\psi}(\varphi,m)U(\varphi,m)d\varphi.
\end{align*}
Letting $n\to\infty$ leads to
\begin{align*}
\int_{0}^{\xi}\int_0^{m}
(A(q_1)(\varphi,\psi)-A(q_2)(\varphi,\psi))\omega(\varphi,\psi)d\varphi d\psi\ge0,
\end{align*}
which completes the proof due to the arbitrariness of $\omega$.
$\hfill\Box$\vskip 4mm

Below we prove the following result for \eqref{strong} similar to the Hopf Lemma.

\begin{lemma}
\label{hopflemma}
Assume that $G\subset\mathbb R^2$ is a circle.
Let $q_1,q_2\in C^2(\overline G)$ with $0<\inf_{G}q_k\le\sup_{G}q_k<c_*\,(k=1,2)$
be sub and super solutions to
\begin{align}
\label{strongq}
\pd{^2{A}(q)}{\varphi^2}+\pd{^2{B}(q)}{\psi^2}=0,\quad(\varphi,\psi)\in G,
\end{align}
respectively.
If $q_1=q_2$ at a point $P_0\in\partial G$ and $q_1<q_2$ on $\overline G\setminus\{P_0\}$,
then $\pd{q_1}{\nu}>\pd{q_2}{\nu}$ at $P_0$,
where $\nu$ is the outer normal to $\partial G$.
\end{lemma}

\Proof
Assume that $G$ is a circle centered at the origin with radius $r$.
Denote $G^*=\{(\varphi,\psi)\in G:\varphi^2+\psi^2>r^2/4\}$, and set
$Q=A(q_1)-A(q_2)$ on $\overline {G^*}$.
Then $Q\in C^2(\overline {G^*})$ is a subsolution to the following linear equation
\begin{align}
\label{acc-1}
\pd{^2Q}{\varphi^2}+\pd{^2}{\psi^2}(b(\varphi,\psi)Q)=0,\quad(\varphi,\psi)\in {G^*},
\end{align}
where
\begin{gather*}
b(\varphi,\psi)=\int_0^1\frac{B'}{A'}\big(A^{-1}(tA(q_1(\varphi,\psi))+(1-t)A(q_2(\varphi,\psi)))\big)dt,\quad
(\varphi,\psi)\in {G^*}.
\end{gather*}
It holds that $b\in C^2(\overline {G^*})$ and
$b_1\le b\le b_2$ in ${G^*}$
with some positive constants $b_1\le b_2$.
Similar to the proof of Proposition \ref{comparison}, one can show that
the comparison principle holds for the problem of \eqref{acc-1} with Dirichlet boundary condition.
Set
$$
\tilde Q(\varphi,\psi)=\mbox{\rm e}^{-\beta r^2}-\mbox{\rm e}^{-\beta (\varphi^2+\psi^2)},\quad
(\varphi,\psi)\in \overline {G^*},
$$
where $\beta$ is a positive constant to be determined.
Direct calculations show that for $(\varphi,\psi)\in {G^*}$,
\begin{align*}
&\pd{^2\tilde Q}{\varphi^2}+\pd{^2}{\psi^2}(b(\varphi,\psi)\tilde Q)
\\
=&-4\beta^2(\varphi^2+b(\varphi,\psi)\psi^2)\mbox{\rm e}^{-\beta (\varphi^2+\psi^2)}
+2\beta(1+b(\varphi,\psi))\mbox{\rm e}^{-\beta (\varphi^2+\psi^2)}
\\
&\qquad+4\beta\pd{b}{\psi}(\varphi,\psi)\psi\mbox{\rm e}^{-\beta (\varphi^2+\psi^2)}
+\pd{^2b}{\psi^2}(\varphi,\psi)\big(\mbox{\rm e}^{-\beta r^2}-\mbox{\rm e}^{-\beta (\varphi^2+\psi^2)}\big)
\\
\le&\Big(-\min\{1,b_1\}r^2\beta^2+2\Big(1+b_2+2\Big\|\pd{b}{\psi}\Big\|_{L^\infty(G^*)}r\Big)\beta
+\Big\|\pd{^2b}{\psi^2}\Big\|_{L^\infty(G^*)}\Big)\mbox{\rm e}^{-\beta (\varphi^2+\psi^2)}.
\end{align*}
Therefore, there exists a suitable constant $\beta>0$ such that
\begin{align*}
\pd{^2\tilde Q}{\varphi^2}+\pd{^2}{\psi^2}(b(\varphi,\psi)\tilde Q)\le0,\quad
(\varphi,\psi)\in {G^*}.
\end{align*}
Since $q_1<q_2$ on $\overline G\setminus\{P_0\}$,
the comparison principle leads to
\begin{align*}
Q(\varphi,\psi)\le\tau\tilde Q(\varphi,\psi),\quad
(\varphi,\psi)\in {G^*},\quad \tau=\frac{1}{\mbox{\rm e}^{-\beta r^2/4}-\mbox{\rm e}^{-\beta r^2}}
\inf_{\partial G^*{\cap}G}(A(q_2)-A(q_1))>0,
\end{align*}
which, together with $q_1=q_2$ at $P_0$, yields that $\pd{q_1}{\nu}>\pd{q_2}{\nu}$ at $P_0$.
$\hfill\Box$\vskip 4mm

\subsection{Existence of solutions}

\begin{proposition}
\label{fixexistence}
Assume that $0<c_e<c_*$, $m>0$, and $\zeta$ and $\xi$ satisfy \eqref{ahatcz}.
There is a solution
$q\in C^\infty([0,\xi]\times[0,m]\setminus\{(\zeta,m),(\xi,m)\})\cap H^1((0,\xi)\times(0,m))
\cap C^{0,1}([0,\xi]\times[0,m]\setminus\{(\zeta,m)\})$
to the problem \eqref{prob1}--\eqref{prob6} such that
\begin{gather}
\label{fixexistence-00}
c_l<q(\varphi,\psi)<c_e,\quad
(\varphi,\psi)\in(0,\xi)\times[0,m)\cup(0,\zeta)\times\{m\},
\\
\label{fixexistence-0-1}
\pd{q}{\varphi}(\varphi,\psi)>0,\quad
(\varphi,\psi)\in[0,\xi]\times(0,m),
\\
\label{fixexistence-0-2}
\pd{q}{\psi}(\varphi,\psi)>0,\quad
(\varphi,\psi)\in(0,\xi)\times(0,m)\cup(\zeta,\xi)\times\{m\}.
\end{gather}
\end{proposition}

\Proof
It is assumed that $\xi>\zeta$, and the proof for the case $\xi=\zeta$ is simpler.
For each positive integer $n$, consider the following problem
\begin{align}
\label{nrprob1}
&\pd{^2{A}(q_n)}{\varphi^2}+\pd{^2{B}(q_n)}{\psi^2}=0,
\quad&&(\varphi,\psi)\in G_n,
\\
\label{nrprob2}
&\pd{{A}(q_n)}{\varphi}(0,\psi)=\frac{1}{R_0q_n(0,\psi)\rho(q^2_n(0,\psi))},
\quad&&\psi\in(0,m),
\\
\label{nrprob3}
&\pd{{B}(q_n)}{\psi}(\varphi,0)=0,\quad&&\varphi\in(0,\xi),
\\
\label{nrprob4}
&\pd{{B}(q_n)}{\psi}(\varphi,m)=0,\quad&&\varphi\in(0,\zeta),
\\
\label{nrprob5}
&q_n(\varphi,f_n(\varphi))=c_e,\quad&&\varphi\in(\zeta,\xi),
\\
\label{nrprob6}
&q_n(\xi,\psi)=c_e,\quad&&\psi\in(0,m-2\delta/n),
\end{align}
where $f_n$, $G_n$ and $\delta$ are defined in the proof of Proposition \ref{linearp}.
Set
\begin{align*}
{\mathscr C}=\Big\{{g}\in C^{1,\alpha}([0,m]):&c_l\le{g}\le c_e\mbox{ and }
{g}'\ge0\mbox{ in }(0,m),g'(0)=g'(m)=0,
\\
&\qquad\|g\|_{C^{\alpha}([0,m])}\le M_1,
\|g\|_{C^{1,\alpha}([0,m])}\le M_2\Big\},
\end{align*}
where $\alpha\in(0,1)$ and $M_1,M_2>0$ are constants to be determined.
For each ${g}\in{\mathscr C}$, consider the problem
\begin{align}
\label{anrprob1}
&\pd{^2{A}({\mathscr Q})}{\varphi^2}+\pd{^2{B}({\mathscr Q})}{\psi^2}=0,
\quad&&(\varphi,\psi)\in G_n,
\\
\label{anrprob2}
&\pd{{A}({\mathscr Q})}{\varphi}(0,\psi)=\frac{1}{R_0{g}(\psi)\rho({g}^2(\psi))},
\quad&&\psi\in(0,m),
\\
\label{anrprob3}
&\pd{{B}({\mathscr Q})}{\psi}(\varphi,0)=0,\quad&&\varphi\in(0,\xi),
\\
\label{anrprob4}
&\pd{{B}({\mathscr Q})}{\psi}(\varphi,m)=0,\quad&&\varphi\in(0,\zeta),
\\
\label{anrprob5}
&{\mathscr Q}(\varphi,f_n(\varphi))=c_e,\quad&&\varphi\in(\zeta,\xi),
\\
\label{anrprob6}
&{\mathscr Q}(\xi,\psi)=c_e,\quad&&\psi\in(0,m-2\delta/n).
\end{align}
It is noted that
$$
\overline {\mathscr Q}(\varphi,\psi)=c_e,\quad
\underline {\mathscr Q}(\varphi,\psi)=A^{-1}\Big(A(c_e)
-\frac{\xi-\varphi}{{R_0}c_l\rho(c_l^2)}\Big),
\quad (\varphi,\psi)\in \overline G_n
$$
are super and sub solutions to the problem \eqref{anrprob1}--\eqref{anrprob6}, respectively.
By a standard fixed point argument and the theory on elliptic
equations (see, e.g., \cite{GT}$\S\,15.5$ and \cite{Lieberman4}),
one can show that the problem \eqref{anrprob1}--\eqref{anrprob6}
admits a unique solution ${\mathscr Q}\in C^\infty(G_n)\cap C^{2}(\overline G_n)$ satisfying
\begin{align}
\label{fixexistence-1}
c_l\le A^{-1}\Big(A(c_e)-\frac{\xi-\varphi}{{R_0}c_l\rho(c_l^2)}\Big)
\le {\mathscr Q}(\varphi,\psi)\le c_e,
\quad(\varphi,\psi)\in\overline G_n.
\end{align}
Thanks to the Harnack inequality, there exist two constants $\tilde\alpha\in(0,1)$ and
$\tilde M_1>0$ depending only on $m$, $\zeta$, $R_0$, $c_e$ and $\gamma$, such that
\begin{align}
\label{fixexistence-2-alpha}
\|{\mathscr Q}\|_{C^{\tilde\alpha}([0,\zeta/2]\times[0,m])}\le \tilde M_1.
\end{align}
The Schauder theory gives
\begin{align}
\label{fixexistence-2}
\pd{\mathscr Q}\psi(0,0)=\pd{\mathscr Q}\psi(0,m)=0,\quad
\|{\mathscr Q}\|_{C^{1,\alpha}([0,\zeta/2]\times[0,m])}\le \tilde M_2,\quad
\|{\mathscr Q}\|_{C^{2,\alpha}([0,\zeta/2]\times[0,m])}\le \tilde M_3,
\end{align}
where $\tilde M_2,\tilde M_3>0$ are constants depending only on $m$, $\zeta$, $R_0$, $c_e$, $\gamma$ and $M_1$,
while $\tilde M_3$ also on $M_2$.
Multiplying \eqref{anrprob1} by ${\mathscr Q}-c_e$ and then integrating over $G_n$ by parts,
one gets from  \eqref{anrprob2}--\eqref{fixexistence-1} that
\begin{align}
\label{fixexistence-2-L2}
\|{\mathscr Q}\|_{H^1(G_n)}\le \tilde M_4,
\end{align}
where $\tilde M_4$ is a positive constant depending only on $m$, $\xi$, $R_0$, $c_e$ and $\gamma$.

Denote ${E}(s)={A}({B}^{-1}(s))$ for $s<0$. Then
${E}'(s)>0$, ${E}''(s)<0$ and ${E}'''(s)<0$ for $s<0$.
Set
$w=\pd{{B}({\mathscr Q})}{\varphi}$ and $z=\pd{{B}({\mathscr Q})}{\psi}$ on $\overline G_n$.
Then $w,z\in C^\infty(G_n)\cap C^{1}(\overline G_n)$ solve
\begin{align*}
&p_{1}(\varphi,\psi)\pd{^2w}{\varphi^2}
+\pd{^2w}{\psi^2}
+p_{2}(\varphi,\psi)\pd{w}{\varphi}
+p_{3}(\varphi,\psi){w}=0,\quad(\varphi,\psi)\in G_n,
\\
&p_{1}(\varphi,\psi)\pd{^2z}{\varphi^2}
+\pd{^2z}{\psi^2}
+p_{4}(\varphi,\psi)\pd{z}{\varphi}
+p_{5}(\varphi,\psi)\pd{z}{\psi}
+p_{6}(\varphi,\psi){z}=0,\quad(\varphi,\psi)\in G_n,
\\
&w(0,\psi)=\frac{B'({\mathscr Q}(0,\psi))}{R_0{g}(\psi)\rho({g}^2(\psi))A'({\mathscr Q}(0,\psi))}>0,\quad&&\hskip-50mm\psi\in(0,m),
\\
&\pd{z}{\varphi}(0,\psi)+p_7(\psi)z(0,\psi)=-\frac{\rho({g}^2(\psi)){g}'(\psi)}{R_0{g}^2(\psi)}\le0,
&&\hskip-50mm\psi\in(0,m),
\\
&\pd{w}{\psi}(\varphi,0)=0,\quad z(\varphi,0)=0,\quad&&\hskip-50mm\varphi\in(0,\xi),
\\
&\pd{w}{\psi}(\varphi,m)=0,\quad z(\varphi,m)=0,\quad&&\hskip-50mm\varphi\in(0,\zeta),
\\
&w(\varphi,f_n(\varphi))\ge0,\quad z(\varphi,f_n(\varphi))\ge0,\quad&&\hskip-50mm\varphi\in(\zeta,\xi),
\\
&w(\xi,\psi)\ge0,\quad z(\xi,\psi)=0,\quad&&\hskip-50mm\psi\in(0,m-2\delta/n),
\end{align*}
where $p_{k}\in L^\infty(G_n)\,(1\le k\le6)$ are defined by
\begin{gather*}
p_{1}={E}'({B}({\mathscr Q}))>0,
\quad
p_{2}=3{E}''({B}({\mathscr Q})){\pd{{B}({\mathscr Q})}{\varphi}},
\quad
p_{3}={E}'''({B}({\mathscr Q}))\Big(\pd{{B}({\mathscr Q})}{\varphi}\Big)^2\le0,
\\
p_{4}=2\frac{E''(B({\mathscr Q}))}{E'(B({\mathscr Q}))}\,\pd{A({\mathscr Q})}{\varphi},\quad
p_{5}=-\frac{E''(B({\mathscr Q}))}{E'(B({\mathscr Q}))}\,{\pd{B({\mathscr Q})}{\psi}},
\\
p_{6}=\Big(\frac{E'''(B({\mathscr Q}))}{(E'(B({\mathscr Q})))^2}-
\frac{(E''(B({\mathscr Q})))^2}{(E'(B({\mathscr Q})))^3}\Big)
\Big(\pd{A({\mathscr Q})}{\varphi}\Big)^2\le0,
\end{gather*}
and
\begin{align*}
p_7(\psi)=\frac{E''(B({\mathscr Q}(0,\psi)))}{R_0{g}(\psi)\rho({g}^2(\psi))(E'(B({\mathscr Q}(0,\psi))))^2}<0,
\quad\psi\in(0,m).
\end{align*}
The comparison principle yields
$w\ge0$ and $z\ge0$ on $\overline G_n$.
Hence
\begin{gather}
\label{nrprob10}
\pd{{\mathscr Q}}{\varphi}(\varphi,\psi)\ge0,
\quad
\pd{{\mathscr Q}}{\psi}(\psi,\psi)\ge0,
\quad
(\varphi,\psi)\in\overline G_n.
\end{gather}

Now take $\alpha=\tilde\alpha$, $M_1=\tilde M_1$ and $M_2=\tilde M_2$.
It follows from \eqref{fixexistence-1}--\eqref{fixexistence-2} and \eqref{nrprob10} that
one can define a mapping $J$ from ${\mathscr C}$ to itself
by $J({g})={\mathscr Q}(0,\cdot)\big|_{[0,m]}$.
It follows from \eqref{fixexistence-2} that $J$ is compact.
One can prove the continuity of $J$ by using its compactness
and the uniqueness result for the problem \eqref{anrprob1}--\eqref{anrprob6}
(see, e.g., \cite{WX2}Proposition 4.7).
The Schauder fixed point theorem yields that
$J$ admits a fixed point. Hence there exists a solution $q_n\in C^\infty(G_n)\cap C^{2}(\overline G_n)$
to the problem \eqref{nrprob1}--\eqref{nrprob6}.
Furthermore, it follows from \eqref{fixexistence-2-L2},
\eqref{fixexistence-1} and \eqref{nrprob10} that $\|q_n\|_{H^1(G_n)}\le\tilde M_4$, and
\begin{gather}
\label{fixexistence-3}
c_l\le q_n(\varphi,\psi)\le c_e,
\quad\pd{q_n}{\varphi}(\varphi,\psi)\ge0,
\quad
\pd{q_n}{\psi}(\psi,\psi)\ge0,
\quad(\varphi,\psi)\in\overline G_n.
\end{gather}
Then one can get by a standard limit process that
the problem \eqref{prob1}--\eqref{prob6} admits a solution $q\in H^{1}((0,\xi)\times(0,m))$ satisfying
\begin{gather}
\label{fixexistence-5}
c_l\le q(\varphi,\psi)\le c_e,
\quad\pd{q}{\varphi}(\varphi,\psi)\ge0,
\quad
\pd{q}{\psi}(\psi,\psi)\ge0,
\quad(\varphi,\psi)\in(0,\xi)\times(0,m).
\end{gather}
The Schauder theory shows that
$q\in C^\infty([0,\xi]\times[0,m]\setminus\{(\zeta,m),(\xi,m)\})\cap C^{0,1}([0,\xi]\times[0,m]\setminus\{(\zeta,m)\})$.
Lemma \ref{hopflemma}, \eqref{prob1}--\eqref{prob6} and \eqref{fixexistence-5} imply \eqref{fixexistence-00} and
\begin{gather*}
\pd{q}{\varphi}(0,\cdot)\Big|_{(0,m)}>0,\quad
\pd{q}{\varphi}(\xi,\cdot)\Big|_{(0,m)}>0,\quad
\pd{q}{\psi}(\cdot,m)\Big|_{(\zeta,\xi)}>0.
\end{gather*}
It suffices to verify \eqref{fixexistence-0-1} and \eqref{fixexistence-0-2} for $(\varphi,\psi)\in(0,\xi)\times(0,m)$.
Set $\tilde w=\pd{{B}(q)}{\varphi}$ and $\tilde z=\pd{{B}(q)}{\psi}$
in $(0,\xi)\times(0,m)$.
Then $\tilde w,\tilde z\in C^\infty([0,\xi]\times[0,m]\setminus\{(\zeta,m),(\xi,m)\})$ solve
\begin{gather*}
\tilde p_{1}(\varphi,\psi)\pd{^2\tilde w}{\varphi^2}
+\pd{^2\tilde w}{\psi^2}
+\tilde p_{2}(\varphi,\psi)\pd{\tilde w}{\varphi}
+\tilde p_{3}(\varphi,\psi){\tilde w}=0,\quad(\varphi,\psi)\in(0,\xi)\times(0,m),
\\
\tilde p_{1}(\varphi,\psi)\pd{^2\tilde z}{\varphi^2}
+\pd{^2\tilde z}{\psi^2}
+\tilde p_{4}(\varphi,\psi)\pd{\tilde z}{\varphi}
+\tilde p_{5}(\varphi,\psi)\pd{\tilde z}{\psi}
+\tilde p_{6}(\varphi,\psi){\tilde z}=0,\quad(\varphi,\psi)\in(0,\xi)\times(0,m),
\end{gather*}
where $\tilde p_{k}\in C^\infty((0,\xi)\times(0,m))\,(1\le k\le6)$ are defined 
similarly as $p_k$.
The strong maximum principle and \eqref{fixexistence-5} yield
$\tilde w>0$ and $\tilde z>0$ in $(0,\xi)\times(0,m)$.
$\hfill\Box$\vskip 4mm

The following proposition shows the optimal H\"older continuity
of solutions to the problem \eqref{prob1}--\eqref{prob6} obtained in Proposition \ref{fixexistence}.

\begin{proposition}
\label{fixexistence2}
Assume that $0<c_e<c_*$, $m>0$, and $\zeta$ and $\xi$ satisfy \eqref{ahatcz}.
Let $q$ be the solution to the problem \eqref{prob1}--\eqref{prob6} obtained in Proposition \ref{fixexistence}.
Then $q\in C^{\alpha}([0,\xi]\times[0,m])$ for each exponent $\alpha\in(0,1/2)$,
and there exists a constant $M>0$ depending only on $m$, $\zeta$, $\xi$, $R_0$, $c_e$, $\gamma$ and $\alpha$
such that $\|q\|_{C^{\alpha}([0,\xi]\times[0,m])}\le M$.
\end{proposition}

\Proof
For $\tau\in(0,\delta]$, set $\Omega_\tau=\{(\varphi,\psi)\in\mathbb R^2:(\varphi-\zeta_n)^2+(\psi-m)^2<\tau^2,\psi<m\}$.
Similar to the proof of Proposition \ref{comparison}, one can prove that the comparison principle
holds for the following problem
\begin{align}
\label{new1}
&\pd{^2(\tilde Q)}{\varphi^2}+\pd{^2B(A^{-1}(\tilde Q))}{\psi^2}=0,
\quad&&(\varphi,\psi)\in \Omega_\tau\cap G_n,
\\
&\pd{\tilde Q}{\psi}(\varphi,m)=0,\quad&&\varphi\in(\zeta_n-\tau,\zeta),
\\
&\tilde Q(\varphi,f_n(\varphi))=A(\tilde c),\quad&&\varphi\in(\zeta,\zeta_n+\tau),
\\
\label{new2}
&\tilde Q(\varphi,\psi)=g(\varphi,\psi),\quad&&(\varphi,\psi)\in\partial\Omega_\tau\cap G_n,
\end{align}
where $n>\delta/\tau$,
$\tilde c\in[c_l,c_e]$ is a constant,
and $g\in C(\overline G_n)$.
Note that $A(q_n)$ is a supersolution to the problem \eqref{new1}--\eqref{new2}
if $A(q_n)\ge g$ on $\partial\Omega_\tau\cap G_n$,
where $q_n\in C^\infty(G_n)\cap C^{2}(\overline G_n)$ is the solution to the problem \eqref{nrprob1}--\eqref{nrprob6}.

We construct subsolutions to the quasilinear problem \eqref{new1}--\eqref{new2} similar to
the ones in the proof of Proposition \ref{linearpp}.
It is noted that $({B'}/{A'})'>0$ in $(0,c_*)$.
For $\tilde Q\in C^2(\overline \Omega_\tau\cap \overline G_n)$ satisfying
$A(c_l/2)\le \tilde Q\le A(c_e)$ in $\Omega_\tau\cap G_n$,
one has
\begin{align}
\label{acc-9}
\pd{^2(\tilde Q)}{\varphi^2}+\pd{^2B(A^{-1}(\tilde Q))}{\psi^2}
=&\pd{^2(\tilde Q)}{\varphi^2}+\frac{B'}{A'}(A^{-1}(\tilde Q))\pd{^2\tilde Q}{\psi^2}
+\frac{({B'}/{A'})'(A^{-1}(\tilde Q))}{A'(A^{-1}(\tilde Q))}\Big(\pd{\tilde Q}{\psi}\Big)^2
\nonumber
\\
\ge&\pd{^2(\tilde Q)}{\varphi^2}+\frac{B'}{A'}(A^{-1}(\tilde Q))\pd{^2\tilde Q}{\psi^2},
\quad(\varphi,\psi)\in \Omega_\tau\cap G_n,
\end{align}
\vskip-6mm
\begin{align}
\label{acc-10}
\frac{B'(c_l/2)}{A'(c_l/2)}\le\frac{B'}{A'}(A^{-1}(\tilde Q))\le\frac{B'(c_e)}{A'(c_e)},\quad(\varphi,\psi)\in \Omega_\tau\cap G_n.
\end{align}
Due to \eqref{acc-10}, as shown in the proof of Proposition \ref{linearpp},
there exist an exponent ${\alpha_0}\in(0,1/2)$
and a constant $\sigma>0$, which depend only on $c_e$ and $\gamma$,
such that for
\begin{align}
\label{new3}
\check Q(\varphi,\psi)=A(\tilde c)+\check M r^{\alpha_0}(\varphi,\psi)\sin{\beta}(\varphi,\psi),
\quad(\varphi,\psi)\in \overline \Omega_\tau\cap \overline G_n
\end{align}
satisfying
\begin{align}
\label{acc-12}
A(c_l/2)\le \check Q\le A(c_e)\mbox{ in }\Omega_\tau\cap G_n
\quad\mbox{and}\quad
{\displaystyle\mbox{osc}}_{\Omega_\tau{\cap} G_n}\frac{B'}{A'}(A^{-1}(\check Q))\le\sigma,
\end{align}
it holds that
\begin{align}
\label{acc-14}
\pd{^2(\check Q)}{\varphi^2}+\frac{B'}{A'}(A^{-1}(\check Q))\pd{^2\check Q}{\psi^2}\ge0,
\quad(\varphi,\psi)\in \Omega_\tau\cap G_n,
\end{align}
where $\check M>0$ is a constant,
$$
r(\varphi,\psi)=
\sqrt{{B'(\tilde c)}(\varphi-\zeta_n)^2/{A'(\tilde c)}+(\psi-m)^2},\quad
(\varphi,\psi)\in[0,\xi]\times[0,m],
$$
$$
\beta(\varphi,\psi)=\left\{
\begin{aligned}
&\frac13\arctan\frac{\sqrt{A'(\tilde c)}(\psi-m)}{\sqrt{B'(\tilde c)}(\varphi-\zeta_n)}-\frac{1}{6}\pi,
&&\quad \mbox{ if }(\varphi,\psi)\in(\zeta_n,\xi]\times[0,m],
\\
&-\frac{1}{3}\pi,
&&\quad \mbox{ if }\varphi=\zeta_n,\,\psi\in[0,m],
\\
&\frac13\arctan\frac{\sqrt{A'(\tilde c)}(\psi-m)}{\sqrt{B'(\tilde c)}(\varphi-\zeta_n)}-\frac{1}{2}\pi,
&&\quad \mbox{ if }(\varphi,\psi)\in[0,\zeta_n)\times[0,m].
\end{aligned}
\right.
$$
It follows from \eqref{acc-9} and \eqref{acc-14} that
\begin{align*}
\pd{^2(\check Q)}{\varphi^2}+\pd{^2B(A^{-1}(\check Q))}{\psi^2}\ge0,
\quad(\varphi,\psi)\in \Omega_\tau\cap G_n.
\end{align*}
Since $({B'}/{A'})'>0$ and $A''<0$ in $(0,c_*)$, one has
\begin{gather*}
\underline r \le\frac{r(\varphi,\psi)}{\sqrt{(\varphi-\zeta_n)^2+(\psi-m)^2}}\le\overline r,
\quad(\varphi,\psi)\in(\varphi,\psi)\in(0,\xi)\times(0,m),
\\
{\displaystyle\mbox{osc}}_{\Omega_\tau{\cap} G_n}\frac{B'}{A'}(A^{-1}(\check Q))
\le\mu_0{\displaystyle\mbox{osc}}_{\Omega_\tau{\cap} G_n}\check Q,
\end{gather*}
where
$$
\underline r=\min\Big\{\frac{B'(c_l/2)}{A'(c_l/2)},1\Big\},\quad
\overline r=\max\Big\{\frac{B'(c_e)}{A'(c_e)},1\Big\},\quad
\mu_0=\frac1{A'(c_e)}\max_{[c_l/2,c_e]}\Big(\frac{B'}{A'}\Big)'.
$$
It is clear that
$\sin\beta\le-1/2$ on $[0,\xi]\times[0,m]$.
Therefore,
if $\tau\le {\check M}^{-1/{\alpha_0}} \tau_0$, then \eqref{acc-12} holds,
where
$\tau_0=\min\big\{\delta,{(c_l/2)^{1/{\alpha_0}}}/{\overline r},{(\sigma/\mu_0)^{1/{\alpha_0}}}/{\overline r}\big\}$.
Summing up, for $\tilde c\in[c_l,c_e]$, $\check M>0$ and
$\tau\le {\check M}^{-1/{\alpha_0}} \tau_0$,
$\check Q$ given by \eqref{new3} is a subsolution to the problem \eqref{new1}--\eqref{new2}
if $\check Q\le g$ on $\partial\Omega_\tau\cap G_n$.

Below, we get a lower barrier function of $q_n$ at $(\zeta,m)$
by using a sequence of subsolutions of the form \eqref{new3}
to the problem \eqref{new1}--\eqref{new2}.
For $n>\delta/{\tau_0}$, set
\begin{align*}
\check Q_0(\varphi,\psi)=A({c_0})+r^{{\alpha_0}}(\varphi,\psi)\sin{\beta}(\varphi,\psi),
\quad(\varphi,\psi)\in \overline \Omega_{\tau_0}\cap \overline G_n,
\end{align*}
where
${c_0}=\min\big\{c_e,A^{-1}\big(A(c_l)+(\underline r \tau_0)^{\alpha_0}/2\big)\big\}$.
Then, the above discussion shows that $\check Q_0$ is a subsolution to
\begin{align}
\label{acc-4}
&\pd{^2(\tilde Q)}{\varphi^2}+\pd{^2B(A^{-1}(\tilde Q))}{\psi^2}=0,
\quad&&(\varphi,\psi)\in \Omega_{\tau_0}\cap G_n,
\\
\label{acc-5}
&\pd{\tilde Q}{\psi}(\varphi,m)=0,\quad&&\varphi\in(\zeta_n-\tau_0,\zeta),
\\
\label{acc-6}
&\tilde Q(\varphi,f_n(\varphi))=A({c_0}),\quad&&\varphi\in(\zeta,\zeta_n+\tau_0),
\\
\label{acc-7}
&\tilde Q(\varphi,\psi)=A(c_l),\quad&&(\varphi,\psi)\in\partial \Omega_{\tau_0}\cap G_n.
\end{align}
It follows from \eqref{nrprob1}, \eqref{nrprob4}, \eqref{nrprob5} and \eqref{fixexistence-3}
that $A(q_n)$ is a supersolution to the problem \eqref{acc-4}--\eqref{acc-7}. Therefore,
\begin{align}
\label{acc-18}
A(q_n)(\varphi,\psi)\ge\check Q_0(\varphi,\psi),
\quad(\varphi,\psi)\in \overline \Omega_{\tau_0}\cap \overline G_n.
\end{align}
Take $\tau_1=2^{-1/{\alpha_0}}\tau_0$. For $n>\delta/{\tau_1}$, set
\begin{align*}
\check Q_1(\varphi,\psi)=A(c_1)+2r^{{\alpha_0}}(\varphi,\psi)\sin{\beta}(\varphi,\psi),
\quad(\varphi,\psi)\in \overline \Omega_{\tau_1}\cap \overline G_n,
\end{align*}
where
$c_1=\min\big\{c_e,A^{-1}\big(A({c_0})+(\underline r \tau_1)^{\alpha_0}/2\big)\big\}$.
Then, $\check Q_1$ is a subsolution to
\begin{align*}
&\pd{^2(\tilde Q)}{\varphi^2}+\pd{^2B(A^{-1}(\tilde Q))}{\psi^2}=0,
\quad&&(\varphi,\psi)\in \Omega_{\tau_1}\cap G_n,
\\
&\pd{\tilde Q}{\psi}(\varphi,m)=0,\quad&&\varphi\in(\zeta_n-\tau_1,\zeta),
\\
&\tilde Q(\varphi,f_n(\varphi))=A(c_1),\quad&&\varphi\in(\zeta,\zeta_n+\tau_1),
\\
&\tilde Q(\varphi,\psi)=\check Q_0(\varphi,\psi),\quad&&(\varphi,\psi)\in\partial \Omega_{\tau_1}\cap G_n,
\end{align*}
while $A(q_n)$ is a supersolution to this problem due to
\eqref{nrprob1}, \eqref{nrprob4}, \eqref{nrprob5} and \eqref{acc-18}. Therefore,
\begin{align}
\label{1acc-18}
A(q_n)(\varphi,\psi)\ge\check Q_1(\varphi,\psi),
\quad(\varphi,\psi)\in \overline \Omega_{\tau_1}\cap \overline G_n.
\end{align}
Take $\tau_2=2^{-2/{\alpha_0}}\tau_0$. For $n>\delta/{\tau_2}$, set
\begin{align*}
\check Q_2(\varphi,\psi)=A(c_2)+2^2r^{{\alpha_0}}(\varphi,\psi)\sin{\beta}(\varphi,\psi),
\quad(\varphi,\psi)\in \overline \Omega_{\tau_2}\cap \overline G_n,
\end{align*}
where
$c_2=\min\big\{c_e,A^{-1}\big(A(c_1)+(\underline r \tau_2)^{\alpha_0}\big)\big\}$.
Then, $\check Q_2$ is a subsolution to
\begin{align*}
&\pd{^2(\tilde Q)}{\varphi^2}+\pd{^2B(A^{-1}(\tilde Q))}{\psi^2}=0,
\quad&&(\varphi,\psi)\in \Omega_{\tau_2}\cap G_n,
\\
&\pd{\tilde Q}{\psi}(\varphi,m)=0,\quad&&\varphi\in(\zeta_n-\tau_2,\zeta),
\\
&\tilde Q(\varphi,f_n(\varphi))=A(c_2),\quad&&\varphi\in(\zeta,\zeta_n+\tau_2),
\\
&\tilde Q(\varphi,\psi)=\check Q_1(\varphi,\psi),\quad&&(\varphi,\psi)\in\partial \Omega_{\tau_2}\cap G_n,
\end{align*}
while $A(q_n)$ is a supersolution to this problem due to
\eqref{nrprob1}, \eqref{nrprob4}, \eqref{nrprob5} and \eqref{1acc-18}. Therefore,
\begin{align*}
A(q_n)(\varphi,\psi)\ge\check Q_2(\varphi,\psi),
\quad(\varphi,\psi)\in \overline \Omega_{\tau_2}\cap \overline G_n.
\end{align*}
Repeating the above process, one gets that for each positive integer $k$ and each $n>\delta/{\tau_k}$,
\begin{align*}
A(q_n)(\varphi,\psi)\ge A(c_k)+2^kr^{{\alpha_0}}(\varphi,\psi)\sin{\beta}(\varphi,\psi),
\quad(\varphi,\psi)\in \overline \Omega_{\tau_k}\cap \overline G_n,\quad \tau_k=2^{-{k}/{\alpha_0}}\tau_0,
\end{align*}
where
$c_k=\min\big\{c_e,A^{-1}\big(A(c_{k-1})+2^{k-2}(\underline r \tau_k)^{\alpha_0}\big)\big\}$.
Note that
$2^{k-2}(\underline r \tau_k)^{\alpha_0}=(\underline r \tau_0)^{\alpha_0}/4$ for $k=1,2,\cdots$.
Therefore, there exists a positive integer $k_0$ depending only on $c_e$ and $\gamma$
such that $c_{k_0}=c_e$. Hence for each $n>\delta/{\tau_{k_0}}$,
\begin{align}
\label{kacc-18}
A(q_n)(\varphi,\psi)\ge A(c_e)+2^{k_0}r^{{\alpha_0}}(\varphi,\psi)\sin{\beta}(\varphi,\psi),
\quad(\varphi,\psi)\in \overline \Omega_{\tau_{k_0}}\cap \overline G_n,\quad \tau_{k_0}=2^{-{k_0}/{\alpha_0}}\tau_0.
\end{align}
According to \eqref{fixexistence-3} and \eqref{kacc-18}, one can prove that
$q\in {C^{\alpha_0}([0,\xi]\times[0,m])}$ and
\begin{align}
\label{kacc-19}
\|q\|_{C^{\alpha_0}([0,\xi]\times[0,m])}\le M_0,
\end{align}
where $q$ is the solution to the problem \eqref{prob1}--\eqref{prob6} obtained in Proposition \ref{fixexistence},
and $M_0>0$ is a constant depending only on $m$, $\zeta$, $\xi$, $R_0$, $c_e$ and $\gamma$.
Set $Q=A(q)$ on $[0,\xi]\times[0,m]$.
As the proof of \eqref{acc-9}, $Q$ is a supersolution to the linear equation
\begin{align*}
\pd{^2(\tilde Q)}{\varphi^2}+b(\varphi,\psi)\pd{^2\tilde Q}{\psi^2}=0,
\quad
b(\varphi,\psi)=\frac{B'}{A'}(q(\varphi,\psi)),\quad(\varphi,\psi)\in(0,\xi)\times(0,m).
\end{align*}
For each exponent $\alpha\in(0,1/2)$, 
it follows from \eqref{fixexistence-00}, \eqref{fixexistence-3}, \eqref{kacc-19} and Proposition \ref{linearp}
that
\begin{align}
\label{kacc-19-00}
A(c_e)-\tilde M((\varphi-\zeta)^2+(\psi-m)^2)^{\alpha/2}\le Q(\varphi,\psi)\le A(c_e),
\quad(\varphi,\psi)\in(0,\xi)\times(0,m),
\end{align}
where $\tilde M>0$ is a constant depending only on $m$, $\zeta$, $\xi$, $R_0$, $c_e$, $\gamma$ and $\alpha$.
The H\"older estimates for elliptic equations (\cite{GT}$\S\,8.10$) 
and \eqref{kacc-19-00} complete the proof of the proposition.
$\hfill\Box$\vskip 4mm

\begin{remark}
The regularity in Proposition \ref{fixexistence2}
is almost optimal.
\end{remark}

\section{Free boundary problem}

In this section, we solve the free boundary problem \eqref{prob1}--\eqref{prob7} in ${\mathscr S}$
for $0<c_e<c_*$, ${R_0}{{\vartheta}} {c_l}\rho(c_l^2)<m<{R_0}{{\vartheta}} c_e\rho(c_e^2)$
and $0<\zeta<R_0 c_l$.
To do so, we need the continuous dependence of solutions
to the problem \eqref{prob1}--\eqref{prob6} together with its well-posedness in $\S\, 3$.

\subsection{Continuous dependence of solutions}

\begin{proposition}
\label{continuousdependence}
Assume that $0<c_e<c_*$ and $m>0$.
For given $\zeta_0$ and $\xi_0$ satisfying $0<\zeta_0<R_0c_l$ and $\zeta_0\le\xi_0\le R_0c_l$,
it holds that
$$
\lim_{\stackrel{\zeta\to\zeta_0,\xi\to\xi_0}{0<\zeta<R_0c_l,\zeta\le\xi\le R_0c_l}}q[\zeta,\xi](\varphi,\psi)
=q[\zeta_0,\xi_0](\varphi,\psi)\mbox{ uniformly for } (\varphi,\psi)\in[0,\xi_0)\times[0,m],
$$
where $q[\zeta,\xi]\in C^2([0,\xi]\times[0,m]\setminus\{(\zeta,m),(\xi,m)\})\cap H^1((0,\xi)\times(0,m))
\cap{\mathscr J}$ is the solution to the fixed boundary problem \eqref{prob1}--\eqref{prob6}.
\end{proposition}

\Proof
It is assumed that $\xi_0>\zeta_0$, and the proof for the case $\xi_0=\zeta_0$ is similar.
For convenience, we use $M_i\,(1\le i\le 6)$ to denote generic constants depending only on
$m$, $\zeta_0$, $\xi_0$, $R_0$, $c_e$ and $\gamma$.
Furthermore, a parenthesis after a generic constant means that
this constant depends also on the variables in the parentheses.

First we prove the continuous dependence on $\xi$.
Fix $(\zeta_0+\xi_0)/2\le\xi_1<\xi_2\le R_0c_l$. Denote $q_k=q[\zeta_0,\xi_k]$ for $k=1,2$.
It follows from Propositions \ref{comparison} and \ref{fixexistence} that
\begin{gather}
\label{acc-21}
c_l\le q_2(\varphi,\psi)\le q_1(\varphi,\psi)\le c_e,
\quad (\varphi,\psi)\in[0,\xi_1]\times[0,m],
\\
\label{acc-22}
q_2(\xi_1,\psi)\ge c_e-M_1(\xi_2-\xi_1)=q_1(\xi_1,\psi)-M_1(\xi_2-\xi_1),\quad \psi\in[0,m].
\end{gather}
Set $Q=A(q_1)-A(q_2)$ on $[0,\xi_1]\times[0,m]$.
Then $0\le Q\in C^2([0,\xi_1]\times[0,m]\setminus\{(\zeta_0,m),(\xi_1,m)\})\cap H^1((0,\xi_1)\times(0,m))
\cap C([0,\xi_1]\times[0,m])$ solves
\begin{align*}
&\pd{^2Q}{\varphi^2}+\pd{^2}{\psi^2}(b(\varphi,\psi)Q)=0,
\quad&&(\varphi,\psi)\in(0,\xi_1)\times(0,m),
\\
&\pd{Q}{\varphi}(0,\psi)=-\frac{h(\psi)}{{R_0}}Q(0,\psi),
\quad
Q(\xi_1,\psi)=A(c_e)-A(q_2)(\xi_1,\psi),\quad&&\psi\in(0,m),
\\
&\pd{Q}{\psi}(\varphi,0)=0,\quad&&\varphi\in(0,\xi_1),
\\
&\pd{Q}{\psi}(\varphi,m)=0,\quad&&\varphi\in(0,\zeta_0),
\\
&Q(\varphi,m)=0,\quad&&\varphi\in(\zeta_0,\xi_1),
\end{align*}
where
\begin{gather*}
b(\varphi,\psi)=\int_0^1\frac{B'}{A'}\big(A^{-1}(tA(q_1(\varphi,\psi))+(1-t)A(q_2(\varphi,\psi)))\big)dt,\quad
(\varphi,\psi)\in(0,\xi_1)\times(0,m),
\\
h(\psi)=\int_0^1\frac1{A^{-1}(t{A(q_1(0,\psi))}+(1-t){A(q_2(0,\psi))})}dt,\quad
\psi\in(0,m).
\end{gather*}
Consider its dual problem
\begin{align}
\label{acc-23-1}
&\pd{^2U}{\varphi^2}
+b(\varphi,\psi)\pd{^2U}{\psi^2}=Q(\varphi,\psi),
\quad&&(\varphi,\psi)\in(0,\xi_1)\times(0,m),
\\
\label{acc-23-2}
&\pd{U}{\varphi}(0,\psi)=-\frac{h(\psi)}{{R_0}}U(0,\psi),
\quad U(\xi_1,\psi)=0,\quad&&\psi\in(0,m),
\\
\label{acc-23-3}
&\pd{U}{\psi}(\varphi,0)=0,\quad&&\varphi\in(0,\xi_1),
\\
\label{acc-23-4}
&\pd{U}{\psi}(\varphi,m)=0,\quad&&\varphi\in(0,\zeta_0),
\\
\label{acc-23-6}
&U(\varphi,m)=0,\quad&&\varphi\in(\zeta_0,\xi_1).
\end{align}
It follows from the proof of Proposition \ref{comparison} that
the problem \eqref{acc-23-1}--\eqref{acc-23-6} admits a unique nonpositive solution
$U\in C^2([0,\xi_1]\times[0,m]\setminus\{(\zeta_0,m),(\xi_1,m)\})\cap H^1((0,\xi_1)\times(0,m))\cap C([0,\xi_1]\times[0,m])$.
The classical theory on elliptic equations (\cite{GT}Theorem 8.33) yields that
\begin{align}
\label{acc-23}
0\le\pd{U}{\varphi}(\xi_1,\psi)\le M_2,\quad\psi\in(0,m).
\end{align}
Multiplying the equation of $Q$ by $U$ and then integrating over $(0,\xi_1)\times(0,m)$ by parts, one gets
from \eqref{acc-21}, \eqref{acc-22} and \eqref{acc-23} that
\begin{align}
\label{acc-25}
\int_{0}^{\xi_1}\int_0^{m}Q^2(\varphi,\psi)d\varphi d\psi
=\int_0^{m}Q(\xi_1,\psi)\pd{U}{\varphi}(\xi_1,\psi) d\varphi d\psi
\le M_3(\xi_2-\xi_1).
\end{align}
For each $\eta\in C^2([0,\xi_1]\times[0,m])$ vanishing near $\{0\}\times[0,m]\cup\{(\zeta_0,m)\}\cup\{\xi_1\}\times[0,m]$,
it follows from \eqref{acc-25} and the classical theory on elliptic equations
(\cite[Theorem 8.12]{GT}) that
$\|\eta Q\|_{H^2((0,\xi_1)\times(0,m))}\le M_1(\eta)(\xi_2-\xi_1)^{1/2}$,
which leads to
\begin{align}
\label{acc-26}
\|\eta Q\|_{L^\infty((0,\xi_1)\times(0,m))}\le M_2(\eta)(\xi_2-\xi_1)^{1/2}.
\end{align}
Set
$$
\Lambda=\Big(\mbox{e}\Big\|\pd{^2b}{\psi^2}\Big\|_{L^\infty((0,\zeta_0/2)\times(0,m))}\Big)^{1/2}
+\frac{\mbox{e}}{R_0 c_l},
\quad
\tilde\zeta=\min\Big\{\frac{\zeta_0}2,\frac1{\Lambda}\Big\}.
$$
Then $0\le Q\in C^\infty([0,\tilde\zeta]\times[0,m])$ solves
\begin{align}
\label{acc-22-1}
&\pd{^2\tilde Q}{\varphi^2}+\pd{^2}{\psi^2}(b(\varphi,\psi)\tilde Q)=0,
\quad&&(\varphi,\psi)\in(0,\tilde\zeta)\times(0,m),
\\
\label{acc-22-2}
&\pd{\tilde Q}{\varphi}(0,\psi)=-\frac{h(\psi)}{{R_0}}\tilde Q(0,\psi),
\quad \tilde Q(\tilde\zeta,\psi)=Q(\tilde\zeta,\psi),\quad&&\psi\in(0,m),
\\
\label{acc-22-6}
&\pd{\tilde Q}{\psi}(\varphi,0)=0,\quad
\pd{\tilde Q}{\psi}(\varphi,m)=0,\quad&&\varphi\in(0,\tilde\zeta).
\end{align}
Similar to the proof of Proposition \ref{comparison},
one can show that the comparison principle holds for the problem \eqref{acc-22-1}--\eqref{acc-22-6}.
Set
$$
\bar Q(\varphi,\psi)=\|Q(\tilde\zeta,\cdot)\|_{L^\infty(0,m)}
\big(\mbox{e}^{\Lambda\tilde\zeta}-\mbox{e}^{\Lambda\varphi}\big),\quad(\varphi,\psi)\in[0,\tilde\zeta]\times[0,m].
$$
Then
\begin{align*}
&\pd{^2\bar Q}{\varphi^2}(\varphi,\psi)+\pd{^2}{\psi^2}(b(\varphi,\psi)\bar Q(\varphi,\psi))
\\
=&\|Q(\tilde\zeta,\cdot)\|_{L^\infty(0,m)}
\Big(-\Lambda^2\mbox{e}^{\Lambda\varphi}+\pd{^2b}{\psi^2}(\varphi,\psi)\big(\mbox{e}^{\Lambda\tilde\zeta}-\mbox{e}^{\Lambda\varphi}\big)\Big)
\le0,\quad(\varphi,\psi)\in(0,\tilde\zeta)\times(0,m),
\end{align*}
\vskip-6mm
\begin{align*}
\pd{\bar Q}{\varphi}(0,\psi)+\frac{h(\psi)}{{R_0}}\bar Q(0,\psi)
=\|Q(\tilde\zeta,\cdot)\|_{L^\infty(0,m)}
\Big(-\Lambda\mbox{e}^{\Lambda\varphi}+\frac{h(\psi)}{{R_0}}\big(\mbox{e}^{\Lambda\tilde\zeta}-\mbox{e}^{\Lambda\varphi}\big)\Big)
\le0,\quad\psi\in(0,m).
\end{align*}
Hence $\bar Q$ is a supersolution to the problem \eqref{acc-22-1}--\eqref{acc-22-6}. Thus
\begin{align*}
Q(\varphi,\psi)\le\bar Q(\varphi,\psi)\le\mbox{e}\|Q(\tilde\zeta,\cdot)\|_{L^\infty(0,m)},
\quad(\varphi,\psi)\in(0,\tilde\zeta)\times(0,m),
\end{align*}
which, together with \eqref{acc-21}, \eqref{acc-26} and Proposition \ref{fixexistence2}, leads to that
\begin{align}
\label{acc-27}
\lim_{\xi_2-\xi_1\to0}\|q_1-q_2\|_{L^\infty((0,\xi_1)\times(0,m))}=0.
\end{align}

Turn to the continuous dependence on $\zeta$.
Fix $\zeta_0/2\le\zeta_1<\zeta_2\le(\zeta_0+\xi_0)/2$.
Denote $\check q_k=q[\zeta_k,\xi_0]$ for $k=1,2$ and $\check Q=A(\check q_1)-A(\check q_2)$.
For each exponent $\alpha\in(0,1/2)$,
$\check Q\in C^2([0,\xi_0]\times[0,m]\setminus\{(\zeta_1,m),(\zeta_2,m),(\xi_0,m)\})\cap H^1((0,\xi_0)\times(0,m))
\cap C^\alpha([0,\xi_0]\times[0,m])$ solves
\begin{align*}
&\pd{^2\check Q}{\varphi^2}+\pd{^2}{\psi^2}(\check b(\varphi,\psi)\check Q)=0,
\quad&&(\varphi,\psi)\in(0,\xi_0)\times(0,m),
\\
&\pd{\check Q}{\varphi}(0,\psi)=-\frac{\check h(\psi)}{{R_0}}\check Q(0,\psi),
\quad \check Q(\xi_0,\psi)=0,\quad&&\psi\in(0,m),
\\
&\pd{\check Q}{\psi}(\varphi,0)=0,\quad&&\varphi\in(0,\xi_0),
\\
&\pd{\check Q}{\psi}(\varphi,m)=0,\quad&&\varphi\in(0,\zeta_1),
\\
&\check Q(\varphi,m)=A(\check q_1)(\varphi,m)-A(\check q_2)(\varphi,m),\quad&&\varphi\in(\zeta_1,\xi_0),
\end{align*}
where
\begin{gather*}
\check b(\varphi,\psi)=\int_0^1\frac{B'}{A'}\big(A^{-1}(tA(\check q_1(\varphi,\psi))+(1-t)A(\check q_2(\varphi,\psi)))\big)dt,\quad
(\varphi,\psi)\in(0,\xi_0)\times(0,m),
\\
\check h(\psi)=\int_0^1\frac1{A^{-1}(t{A(\check q_1(0,\psi))}+(1-t){A(\check q_2(0,\psi))})}dt,\quad
\psi\in(0,m).
\end{gather*}
Consider its dual problem
\begin{align*}
&\pd{^2\check U}{\varphi^2}
+\check b(\varphi,\psi)\pd{^2\check U}{\psi^2}=\check Q(\varphi,\psi),
\quad&&(\varphi,\psi)\in(0,\xi_0)\times(0,m),
\\
&\pd{\check U}{\varphi}(0,\psi)=-\frac{\check h(\psi)}{{R_0}}\check U(0,\psi),
\quad \check U(\xi_0,\psi)=0,\quad&&\psi\in(0,m),
\\
&\pd{\check U}{\psi}(\varphi,0)=0,\quad&&\varphi\in(0,\xi_0),
\\
&\pd{\check U}{\psi}(\varphi,m)=0,\quad&&\varphi\in(0,\zeta_1),
\\
&\check U(\varphi,m)=0,\quad&&\varphi\in(\zeta_1,\xi_0),
\end{align*}
which admits a unique solution
$\check U\in C^2([0,\xi_0]\times[0,m]\setminus\{(\zeta_1,m),(\xi_0,m)\})\cap H^1((0,\xi_0)\times(0,m))
\cap C^\alpha([0,\xi_0]\times[0,m])$ from the proof of Proposition \ref{comparison}.
For $(\varphi,\psi)\in(0,\xi_0)\times(0,m)$ satisfying $(\varphi-(\zeta_1+\zeta_2)/2)^2+(\psi-m)^2\le(\zeta_2-\zeta_1)^2$,
it follows from the H\"older estimates for elliptic equations
(\cite{GT}$\S\,8.11$) that
\begin{gather*}
|\nabla\check Q(\varphi,\psi)|\le M_4(\alpha)\big((\varphi-\zeta_2)^2+(\psi-m)^2\big)^{-(1-\alpha)/2},
\\
|\nabla\check U(\varphi,\psi)|\le M_5(\alpha)\big((\varphi-\zeta_1)^2+(\psi-m)^2\big)^{-(1-\alpha)/2}.
\end{gather*}
Choose $0\le\omega\in C^1([0,\xi_0]\times[0,m])$ such that for $(\varphi,\psi)\in(0,\xi_0)\times(0,m)$,
$$
\omega(\varphi,\psi)=
\left\{
\begin{aligned}
&1,\quad\mbox{if }(\varphi-(\zeta_1+\zeta_2)/2)^2+(\psi-m)^2\ge4(\zeta_2-\zeta_1)^2,
\\
&0,\quad\mbox{if }(\varphi-(\zeta_1+\zeta_2)/2)^2+(\psi-m)^2\le(\zeta_2-\zeta_1)^2,
\end{aligned}
\right.
$$
\vskip-2mm
\begin{gather*}
|\nabla\omega(\varphi,\psi)|\le\frac{4}{\zeta_2-\zeta_1},
\quad
\Big|\pd{^2\omega}{\varphi^2}(\varphi,\psi)\Big|+
\Big|\pd{^2\omega}{\psi^2}(\varphi,\psi)\Big|\le\frac{8}{(\zeta_2-\zeta_1)^2}.
\end{gather*}
Multiplying the equation of $\check Q$ by $\omega\check U$ and then integrating over $(0,\xi_0)\times(0,m)$ by parts, one gets
that
\begin{align*}
\int_{0}^{\xi_0}\int_0^{m}\omega(\varphi,\psi)\check  Q^2(\varphi,\psi)d\varphi d\psi
\le M_6(\alpha)(\zeta_2-\zeta_1)^{2\alpha}.
\end{align*}
A similar argument as \eqref{acc-27} leads to
$\lim_{\zeta_2-\zeta_1\to0}\|\check q_1-\check q_2\|_{L^\infty((0,\xi_0)\times(0,m))}=0$.
$\hfill\Box$\vskip 4mm

\subsection{Solvability and properties of solutions}

Using Propositions \ref{comparison}, \ref{fixexistence}, \ref{fixexistence2}, \ref{continuousdependence} and Lemma \ref{hopflemma},
one can prove the following three lemmas.

\begin{lemma}
\label{uniqueness}
Assume that $0<c_e<c_*$ and ${R_0}{{\vartheta}} {c_l}\rho(c_l^2)<m<{R_0}{{\vartheta}} c_e\rho(c_e^2)$.

{\rm(i)}
For $0<\zeta<R_0c_l$,
there is at most one solution to the problem \eqref{prob1}--\eqref{prob7} in ${\mathscr S}$.

{\rm(ii)}
For $\hat\zeta<\zeta<R_0c_l$, there is not a solution to the problem \eqref{prob1}--\eqref{prob7}
in ${\mathscr S}$.
\end{lemma}

\Proof
For $0<\zeta<R_0c_l$, assume that $(q_1,\xi_1),(q_2,\xi_2)\in{\mathscr S}$ are two solutions
to the problem \eqref{prob1}--\eqref{prob7}.
We first prove $\xi_1=\xi_2$ by contradiction. If not,
it is assumed that $\xi_1<\xi_2$ without loss of generality.
Then, Propositions  \ref{comparison}, \ref{fixexistence}, \ref{fixexistence2} lead to
$q_2(\xi_1,\psi)<c_e=q_1(\xi_1,\psi)$ for $\psi\in(0,m)$.
It follows from this estimate, Proposition \ref{comparison} and Lemma \ref{hopflemma} that
$q_2(0,\psi)<q_1(0,\psi)$ for $\psi\in[0,m]$,
which contradicts \eqref{prob7}. Hence $\xi_1=\xi_2$. And Proposition \ref{comparison} shows $q_1=q_2$.

For $\hat\zeta<\zeta<R_0c_l$,
assume that $(q,\xi)\in{\mathscr S}$ is a solution to the problem \eqref{prob1}--\eqref{prob7}.
Then $\xi\ge\zeta>\hat\zeta$. Propositions  \ref{comparison}, \ref{fixexistence}, \ref{fixexistence2}
and Lemma \ref{symmetricexistence} lead to
$q(\hat\zeta,\psi)<c_e=\hat q(\hat\zeta)$ for $\psi\in(0,m)$,
which, together with Proposition \ref{comparison} and Lemma \ref{hopflemma}, shows that
$q(0,\psi)<\hat q(0)$ for $\psi\in[0,m]$.
Hence
\begin{align*}
\int_{0}^{m}\frac1{q(0,\psi)\rho(q^2(0,\psi))}d\psi>\frac{m}{\hat q(0)\rho(\hat q^2(0))}={R_0}{{\vartheta}},
\end{align*}
which contradicts \eqref{prob7}.
$\hfill\Box$\vskip 4mm

\begin{lemma}
\label{existencelemma1}
Assume that for $\zeta_0\in(0,\hat\zeta]$,
the problem \eqref{prob1}--\eqref{prob7} with $\zeta=\zeta_0$ admits a solution $(q_0,\xi_0)\in{\mathscr S}$
with $\xi_0<R_0c_l$.
Then, there exists $\tau\in(0,\zeta_0)$, such that for each
$\zeta\in(\zeta_0-\tau,\zeta_0)$,
the problem \eqref{prob1}--\eqref{prob7} admits a solution
$(q,\xi)\in{\mathscr S}$.
\end{lemma}

\Proof
For $\zeta=\zeta_0$,
denote $q_1,q_2\in{\mathscr J}$ to be the solutions to the problem \eqref{prob1}--\eqref{prob6} with
$\xi=R_0c_l$ and $\xi=(\xi_0+R_0c_l)/2$, respectively. Then,
Propositions  \ref{comparison}, \ref{fixexistence}, \ref{fixexistence2} and Lemma \ref{hopflemma} yield
$c_l\le q_1(0,\psi)<q_2(0,\psi)<q_0(0,\psi)$ for $\psi\in[0,m]$.
Due to Propositions  \ref{comparison}, \ref{fixexistence}, \ref{fixexistence2}
and \ref{continuousdependence}, there exists $\tau\in(0,\zeta_0)$, such that for
$\zeta\in(\zeta_0-\tau,\zeta_0)$, the solution $q_{+,\zeta}\in{\mathscr J}$
to the problem \eqref{prob1}--\eqref{prob6} with $\xi=(\xi_0+R_0c_l)/2$
satisfies
$c_l\le q_{+,\zeta}(0,\psi)<q_0(0,\psi)$ for $\psi\in[0,m]$.
Hence
\begin{align}
\label{free1}
\int_{0}^{m}\frac1{q_{+,\zeta}(0,\psi)\rho(q_{+,\zeta}^2(0,\psi))}d\psi>\int_{0}^{m}\frac1{q_0(0,\psi)\rho(q^2_0(0,\psi))}d\psi={R_0}{{\vartheta}}.
\end{align}
For $\zeta\in(\zeta_0-\tau,\zeta_0)$, Propositions  \ref{comparison}, \ref{fixexistence}, \ref{fixexistence2}
and Lemma \ref{hopflemma} show that
the solution $q_{-,\zeta}\in{\mathscr J}$ to the problem \eqref{prob1}--\eqref{prob6} with $\xi=\xi_0$ satisfies
$q_{-,\zeta}(0,\psi)>q_0(0,\psi)$ for $\psi\in[0,m]$,
which yields
\begin{align}
\label{free2}
\int_{0}^{m}\frac1{q_{-,\zeta}(0,\psi)\rho(q_{-,\zeta}^2(0,\psi))}d\psi<\int_{0}^{m}\frac1{q_0(0,\psi)\rho(q^2_0(0,\psi))}d\psi={R_0}{{\vartheta}}.
\end{align}
For $\zeta\in(\zeta_0-\tau,\zeta_0)$, it follows from \eqref{free1}, \eqref{free2},
Propositions \ref{fixexistence}, \ref{fixexistence2} and \ref{continuousdependence}
that there exists $\xi\in(\xi_0,R_0c_l)$ such that
the problem \eqref{prob1}--\eqref{prob7} admits a solution $(q,\xi)\in{\mathscr S}$.
$\hfill\Box$\vskip 4mm

\begin{lemma}
\label{existencelemma2}
Assume that for $\zeta_0\in(0,\hat\zeta)$,
the problem \eqref{prob1}--\eqref{prob7} with $\zeta=\zeta_0$ admits a solution $(q_0,\xi_0)\in{\mathscr S}$.
Then, for each $\zeta\in(\zeta_0,\hat\zeta)$,
the problem \eqref{prob1}--\eqref{prob7} admits a solution $(q,\xi)\in{\mathscr S}$.
\end{lemma}

\Proof
Give $\zeta\in(\zeta_0,\hat\zeta)$.
Denote $q_1,q_2\in{\mathscr J}$ to be the solutions to the problem \eqref{prob1}--\eqref{prob6} with $\xi=\xi_0$
and $\xi=\hat\zeta$, respectively.
Then, Propositions  \ref{comparison}, \ref{fixexistence}, \ref{fixexistence2}
and Lemmas \ref{symmetricexistence}, \ref{hopflemma} yield
$c_l\le q_1(0,\psi)<q_0(0,\psi)$ and $q_2(0,\psi)>\hat q(0)>c_l$ for $\psi\in[0,m]$.
Hence
\begin{gather*}
\int_{0}^{m}\frac1{q_1(0,\psi)\rho(q^2_1(0,\psi))}d\psi>\int_{0}^{m}\frac1{q_0(0,\psi)\rho(q^2_0(0,\psi))}d\psi={R_0}{{\vartheta}},
\\
\int_{0}^{m}\frac1{q_2(0,\psi)\rho(q^2_2(0,\psi))}d\psi<\frac{m}{\hat q(0)\rho(\hat q^2(0))}={R_0}{{\vartheta}}.
\end{gather*}
Therefore, Propositions \ref{fixexistence}, \ref{fixexistence2} and \ref{continuousdependence}
show that
there exists $\xi\in(\hat\zeta,\xi_0)$ such that
the problem \eqref{prob1}--\eqref{prob7} admits a solution $(q,\xi)\in{\mathscr S}$.
$\hfill\Box$\vskip 4mm

We are ready to prove the following existence and nonexistence results for the free boundary problem \eqref{prob1}--\eqref{prob7}
in ${\mathscr S}$.

\begin{theorem}
\label{existence}
For $0<c_e<c_*$ and ${R_0}{{\vartheta}} {c_l}\rho(c_l^2)<m<{R_0}{{\vartheta}} c_e\rho(c_e^2)$,
there exists a constant $\zeta_*\in[0,\hat\zeta)$, depending only on $R_0$, ${\vartheta}$, $c_e$ and $\gamma$,
such that the problem \eqref{prob1}--\eqref{prob7} admits a unique solution $(q[\zeta],\xi[\zeta])\in
C^2([0,\xi]\times[0,m]\setminus\{(\zeta,m),(\xi,m)\})
\cap H^1((0,\xi)\times(0,m))\cap{\mathscr S}$
if $\zeta\in[\zeta_*,\hat\zeta]\cap(0,\hat\zeta]$,
while there is not such a solution 
if $\zeta\in(0,\zeta_*)\cup(\hat\zeta,R_0c_l)$,
where ${c_l}$ and $\hat\zeta$
are given in Theorem \ref{existencephysical} and Lemma \ref{symmetricexistence}, respectively.
If $\zeta_*>0$ additionally, then
$\xi[\zeta_*]=R_0c_l$.
\end{theorem}

\Proof
Set $\zeta_*=\inf{\mathscr E}$, where
\begin{align*}
{\mathscr E}=\big\{\zeta\in(0,R_0c_l]:
\mbox{the problem \eqref{prob1}--\eqref{prob7} admits a unique solution }(q[\zeta],\xi[\zeta])\in{\mathscr S}\big\}.
\end{align*}
It follows from Lemmas \ref{uniqueness},
\ref{existencelemma1}, \ref{existencelemma2} and \ref{symmetricexistence} that
$$
\zeta_*\in[0,\hat\zeta),\quad (\zeta_*,\hat\zeta]\subset{\mathscr E},
\quad
(0,\zeta_*)\cap{\mathscr E}=\emptyset,\quad (\hat\zeta,R_0c_l)\cap{\mathscr E}=\emptyset.
$$
If $\zeta_*>0$ additionally,
Proposition \ref{continuousdependence} and Lemma \ref{existencelemma1}
show that $\zeta_*\in{\mathscr E}$
and $\xi[\zeta_*]=R_0c_l$.
$\hfill\Box$\vskip 4mm

Solutions to the problem \eqref{prob1}--\eqref{prob7} have the following properties.

\begin{theorem}
\label{existencep}
Assume that $0<c_e<c_*$ and ${R_0}{{\vartheta}} {c_l}\rho(c_l^2)<m<{R_0}{{\vartheta}} c_e\rho(c_e^2)$.
For $\zeta\in[\zeta_*,\hat\zeta]\cap(0,\hat\zeta]$,
let $(q[\zeta],\xi[\zeta])\in{\mathscr S}({c_l})$ be the unique solution to the problem \eqref{prob1}--\eqref{prob7},
where $\zeta_*$, ${c_l}$ and $\hat\zeta$
are given in Theorem \ref{existence}, Theorem \ref{existencephysical} and Lemma \ref{symmetricexistence}, respectively.

{\rm(i)}
$q[\zeta],\theta[\zeta]\in C^\infty([0,\xi]\times[0,m]\setminus\{(\zeta,m),(\xi,m)\})\cap H^1((0,\xi)\times(0,m))
\cap C^{0,1}([0,\xi]\times[0,m]\setminus\{(\zeta,m)\})\cap C^{\alpha}([0,\xi]\times[0,m])$
for each exponent $\alpha\in(0,1/2)$, and
\begin{gather*}
c_l<q[\zeta](\varphi,\psi)< c_e,\quad
(\varphi,\psi)\in(0,\xi)\times[0,m)\cup(0,\zeta)\times\{m\},
\\
-{\vartheta}<\theta[\zeta](\varphi,\psi)<0,\quad
(\varphi,\psi)\in[0,\xi]\times(0,m)\cup(\zeta,\xi]\times\{m\},
\\
\pd{q[\zeta]}{\varphi}(\varphi,\psi)>0,\quad
\pd{\theta[\zeta]}{\psi}(\varphi,\psi)<0,\quad
(\varphi,\psi)\in[0,\xi]\times(0,m),
\\
\pd{q[\zeta]}{\psi}(\varphi,\psi)>0,\quad
\pd{\theta[\zeta]}{\varphi}(\varphi,\psi)>0,\quad
(\varphi,\psi)\in(0,\xi)\times(0,m)\cup(\zeta,\xi)\times\{m\},
\end{gather*}
where $\theta[\zeta]$ is the flow angle of $q[\zeta]$.

{\rm(ii)}
$(q[\hat\zeta],\xi[\hat\zeta])=(\hat q,\hat\zeta)$ with $\hat q$ given in Lemma \ref{symmetricexistence}.

{\rm(iii)}
If $\zeta_*>0$ additionally, then $\xi[\zeta_*]=R_0{c_l}$.

{\rm(iv)}
For each $\zeta_0\in[\zeta_*,\hat\zeta]\cap(0,\hat\zeta]$,
\begin{gather*}
\lim_{\stackrel{\zeta\to\zeta_0}{\zeta\in[\zeta_*,\hat\zeta]{\cap}(0,\hat\zeta]}}\xi[\zeta]
=\xi[\zeta_0],
\\
\lim_{\stackrel{\zeta\to\zeta_0}{\zeta\in[\zeta_*,\hat\zeta]{\cap}(0,\hat\zeta]}}(q[\zeta],\theta[\zeta])(\varphi,\psi)
=(q[\zeta_0],\theta[\zeta_0])(\varphi,\psi)\mbox{ uniformly for } (\varphi,\psi)\in[0,\xi[\zeta_0])\times[0,m].
\end{gather*}

{\rm(v)}
For $\zeta_1,\zeta_2\in[\zeta_*,\hat\zeta]\cap(0,\hat\zeta]$ with
$\zeta_1<\zeta_2$, it holds
\begin{gather}
\label{existence-comarison1-0}
\xi[\zeta_1]>\xi[\zeta_2],\quad L[\zeta_1]>L[\zeta_2],
\\
\label{existence-comarison1-1}
q[\zeta_1](\varphi,m)>q[\zeta_2](\varphi,m),\quad\varphi\in[0,\zeta_1),
\\
\label{existence-comarison1-2}
q[\zeta_1](\varphi,0)<q[\zeta_2](\varphi,0),\quad\varphi\in[0,\xi[\zeta_2]],
\end{gather}
where $L[\zeta]$ is the length of the upper wall of the nozzle in the physical plane for the flow $(q[\zeta],\xi[\zeta])$.
\end{theorem}

\Proof
Theorem \ref{existence}, Propositions  \ref{comparison}, \ref{fixexistence}, \ref{fixexistence2} and \ref{continuousdependence},
Lemmas \ref{symmetricexistence} and \ref{hopflemma}
yield (i)--(iv) directly,
and it suffices to verify (v).
Give $\zeta_1,\zeta_2\in[\zeta_*,\hat\zeta]\cap(0,\hat\zeta]$ with $\zeta_1<\zeta_2$.
For convenience, denote $q_k=q[\zeta_k]$ and $\xi_k=\xi[\zeta_k]$ for $k=1,2$.
If $\xi_1\le\xi_2$, then (i) leads to
$q_2(\xi_1,\psi)<c_e=q_1(\xi_1,\psi)$ for $\psi\in(0,m)$.
This estimate, Proposition \ref{comparison} and Lemma \ref{hopflemma} show that
$q_2(0,\psi)<q_1(0,\psi)$ for $\psi\in[0,m]$,
which contradicts \eqref{prob7}. Hence $\xi_1>\xi_2$,
which yields the first estimate in \eqref{existence-comarison1-0}.
It follows from $\xi_1>\xi_2$ and (i) that
\begin{align}
\label{free4}
q_1(\xi_2,\psi)<c_e=q_2(\xi_2,\psi),\quad
\psi\in(0,m).
\end{align}
Set $\Gamma_+=[\zeta_1,\zeta_2)\times\{m\}$, $\Gamma_-=\{\xi_2\}\times[0,m]$, and
$G_\pm=\big\{(\varphi,\psi)\in[0,\xi_2]\times[0,m]:\pm q_1(\varphi,\psi)>\pm q_2(\varphi,\psi)\big\}$.
Then, (i) and \eqref{free4} show that
$G\pm$ are relatively open sets in $[0,\xi_2]\times[0,m]$ and $\Gamma_\pm\subset G_\pm$.
We claim that $G_\pm$ are connected. If not,
let $G_\pm^*$ be a maximal subdomain of $G_\pm$ such that $\Gamma_\pm\cap G_\pm^*=\emptyset$.
Note that $q_1$ and $q_2$ satisfy the same boundary conditions on
$\partial((0,\xi_2)\times(0,m))\setminus(\Gamma_+\cup\Gamma_-)$,
$q_1=q_2$ on $\partial G_\pm^*\cap (0,\xi_2)\times(0,m)$,
and $G_\pm^*$ satisfy the interior cone condition.
Then, the comparison principle, which can be proved similar to Proposition \ref{comparison}, leads to
$\pm q_1\le\pm q_2$ in $G_\pm^*$,
which contradicts $G_\pm^*\subset G_\pm$. Hence
\begin{align}
\label{free5}
G\pm \mbox{ are connected and relatively open sets in } [0,\xi_2]\times[0,m].
\end{align}
Next we show that
\begin{align}
\label{free6}
[0,\xi_2]\times\{0\}\cap\overline G_+=\emptyset
\end{align}
by contradiction.
If not, from \eqref{free5} and $\Gamma_+\subset G_+$,
there exists a $C^1$ simple curve ${\mathscr L}_+^*\subset (0,\xi_2)\times[0,m]\cup\{(0,0),(\xi_2,0)\}$ from $[0,\xi_2]\times\{0\}$ to $\Gamma_+$
such that $q_1\ge q_2$ on ${\mathscr L}_+^*$
and ${\mathscr G}_+$ satisfies the interior cone condition, where
${\mathscr G}_+=\big\{(\varphi,\psi)\in(0,\xi_2)\times(0,m):
(\varphi,\psi) \mbox{ is located at the left side of } {\mathscr L}_+^*\big\}$.
Using the comparison principle, which can be proved in a similar way as Proposition \ref{comparison},
one gets that $q_1\ge q_2$ on $\overline{\mathscr G}_+$.
Then, Lemma \ref{hopflemma} leads to
$q_1(0,\psi)>q_2(0,\psi)$ for $\psi\in(0,m)$,
which contradicts \eqref{prob7}.
Similarly, one can prove that
\begin{align}
\label{free7}
[0,\zeta_1]\times\{m\}\cap\overline G_-=\emptyset.
\end{align}
If not, from \eqref{free5} and $\Gamma_-\subset G_-$,
there exists a $C^1$ simple curve ${\mathscr L}_-^*\subset (0,\xi_2]\times(0,m)\cup [0,\zeta_1]\times\{m\}$ from $\Gamma_-$
to $[0,\zeta_1]\times\{m\}$
such that $q_1\le q_2$ on ${\mathscr L}_-^*$
and ${\mathscr G}_-$ satisfies the interior cone condition, where
${\mathscr G}_-=\big\{(\varphi,\psi)\in(0,\xi_2)\times(0,m):
(\varphi,\psi) \mbox{ is located at the left side of } {\mathscr L}_-^*\big\}$.
Using the comparison principle and Lemma \ref{hopflemma},
one gets that
$q_1(0,\psi)<q_2(0,\psi)$ for $\psi\in(0,m)$,
which contradicts \eqref{prob7}.
It follows from \eqref{free5}--\eqref{free7} that
there exists a suitable constant $\tau\in(0,m)$ such that
\begin{gather*}
q_1\le q_2\mbox{ in } [0,\xi_2]\times[0,\tau]
\quad\mbox{and}\quad
q_1\ge q_2\mbox{ in } [0,\zeta_1]\times[m-\tau,m],
\end{gather*}
which, together with Lemma \ref{hopflemma}, leads to \eqref{existence-comarison1-1} and \eqref{existence-comarison1-2}.
Note that
$$
{L}[\zeta]=\int_{0}^{\zeta}\frac1{q[\zeta](\varphi,m)}d\varphi.
$$
So, the second estimate in \eqref{existence-comarison1-0} follows from the first estimate in \eqref{existence-comarison1-0} and
\eqref{existence-comarison1-1}.
$\hfill\Box$\vskip 4mm

Transforming Theorems \ref{existence} and \ref{existencep} into the physical plane,
one can get Theorem \ref{existencephysical} and \ref{existencephysical2}.

\end{document}